\documentclass[pdf, relax, 11pt]{article}

\usepackage{amsmath,amsfonts,amssymb}

\usepackage{array}
\usepackage{bm}
\usepackage{caption2}
\usepackage{color}
\usepackage{cite}
\usepackage{epic}
\usepackage{footnote}
\usepackage[bottom]{footmisc}
\usepackage{graphics}
\usepackage{graphicx}
\usepackage{import}
\usepackage[utf8]{inputenc}
\usepackage{marginnote}
\usepackage{mathtools}
\usepackage{multicol}  
\usepackage{multirow}
\usepackage{psfrag}
\usepackage{setspace}
\usepackage{tabularx}
\usepackage{url}
\usepackage{algorithmic}

\definecolor{dgreen}{rgb}{0.133333, 0.545098, 0.133333}

\title{A Dynamically Adaptive Sparse Grid Method for Quasi-Optimal Interpolation of Multidimensional Analytic Functions}
\author{Miroslav~K.~Stoyanov, Clayton~G.~Webster}

\newtheorem{theorem}{Theorem}
\newtheorem{algorithm}{Algorithm}
\newtheorem{assume}{Assumption}
\newtheorem{rem}{Remark}
\newtheorem{cor}{Corollary}

\author{Miroslav ~Stoyanov\footnote{
Department of Computational and Applied Mathematics, 
Oak Ridge National Laboratory, One Bethel Valley Road, 
P.O.~Box 2008, MS-6367, Oak Ridge, TN 37831-6164 
({\tt stoyanovmk@ornl.gov})}
\quad
Clayton ~Webster\footnote{
Department of Computational and Applied Mathematics,
Oak Ridge National Laboratory, One Bethel Valley Road,
P.O.~Box 2008, MS-6164, Oak Ridge, TN 37831-6164
({\tt webstercg@ornl.gov})}
}

\newcommand{\T}{\Theta}

\begin{document}
\maketitle

\abstract{
In this work we develop a dynamically adaptive sparse grids (SG) method for quasi-optimal interpolation of multidimensional analytic functions defined over a product of one dimensional bounded domains. The goal of such approach is to construct an interpolant in space that corresponds to the ``best $M$-terms'' based on sharp a priori estimate of polynomial coefficients. In the past, SG methods have been successful in achieving this, with a traditional construction that relies on the solution to a Knapsack problem: only the most profitable hierarchical surpluses are added to the SG. However, this approach requires additional sharp estimates related to the size of the analytic region and the norm of the interpolation operator, i.e., the Lebesgue constant. Instead, we present an iterative SG procedure that adaptively refines an estimate of the region and accounts for the effects of the Lebesgue constant. Our approach does not require any a priori knowledge of the analyticity or operator norm, is easily generalized to both affine and non-affine analytic functions, and can be applied to sparse grids build from one dimensional rules with arbitrary growth of the number of nodes. In several numerical examples, we utilize our dynamically adaptive SG to interpolate quantities of interest related to the solutions of parametrized elliptic and hyperbolic PDEs, and compare the performance of our quasi-optimal interpolant to several alternative SG schemes.}

\section{Introduction}

This paper considers constructing an approximations to multidimensional analytic functions defined over a product of one dimensional bounded domains. The main challenge facing all methods in this context is the \emph{curse of dimensionality}, i.e., the computational complexity of approximation techniques increases exponentially with the number of dimensions. To alleviate the curse, methods have been proposed that reduce the dimensionality of the problem\cite{stoyanov2015gradient,Loeve1987}, reduce the complexity the target function\cite{binev2011convergence,devore2013greedy}, or approximate the function in an optimal polynomial subspace\cite{nobile2008sparse,nobile2008anisotropic,cohen2011analytic}. We take the latter approach and we build upon the recent results in best $M$-terms approximation\cite{cohen2011analytic,TranQuasiOptimal,chkifa2014high}, where the function is projected onto the polynomial space associated with the dominant coefficients of either a Taylor or Legendre expansion. In implementation, finding the optimal space is intractable and instead sharp a priori estimates of the expansion coefficients are used to select a quasi-optimal space. Such approach can achieve sub-exponential convergence rate in the context of both projection, e.g., \cite{Beck2014732,todor2007convergence,DexterCounting} and interpolation, e.g., \cite{chkifa2014high,NobileQuasiOptimalInterpolation}, however, the quasi-optimal methods rely heavily on a priori estimates of the size of the region of analyticity of the function and sharp estimates are available only in few special cases.

Given a suitable polynomial space, orthogonal projection results in the best $L^2$ approximation, however, the projection approach often times comes at a heavy computational cost\cite{Beck2014732,todor2007convergence,DexterCounting}. In contrast, sampling based techniques require only the values of the function at a set of nodes, e.g., Monte Carlo random sampling for computing the statistical moments of a function\cite{niederreiter1992quasi,Fishman_96}, and Sparse Grids (SG) method for high order polynomial approximation\cite{nobile2008sparse,nobile2008anisotropic,gunzburger2014stochastic}, which is the focus on this work. SG sampling does not result in best approximation in the associated polynomial space and the error is magnified by the norm of the SG operator, a.k.a., the Lebesgue constant. However, sampling tends to be computationally cheaper than projection as well as more susceptible to parallelization which usually offsets the moderate increase of the error. In addition, sampling procedures can wrap around simulation software that computes single realization of the function, which simplifies the implementation and allows the use of legacy and third party code.

Sparse grids algorithms construct multidimensional function approximation from a linear combination of tensors of one dimensional interpolation rules. Quasi-optimal SG are traditionally constructed as the solution to a Knapsack problem\cite{bungartz2004sparse,NobileQuasiOptimalInterpolation}, where the selected set of tensors is associated with the largest profit index that is derived from an a priori estimate of the hierarchical surplus, the Lebesgue constant, and the number of samples in a tensor. In the case when the one dimensional rules grow by one node at a time, a near optimal greedy procedure using the Taylor coefficients of the function can construct a suitable approximation\cite{chkifa2013sparse,chkifa2014high}, however, without a priori assumptions, selecting the optimal set of coefficients comes at a very high computational cost.

In this work, we present an iterative procedure for constructing a sequence of SG interpolants with increasing number of nodes and  accuracy, that does not require a priori estimates of the region of analyticity. We focus our attention to the nested SG case, where all nodes associated with one grid are also utilized by the next grid in the sequence, thus reusing all available samples. We review popular one dimensional nested rules such as Clenshaw-Curtis\cite{clenshaw1960method} and Leja\cite{de2004leja,chkifa2013lebesgue} and we present several new rules based on greedy minimization of operator norms. In addition, for any chosen rule and any arbitrary lower (i.e., admissible\cite{bungartz2004sparse}) polynomial space, we present a strategy for selecting the minimal set of tensors that yields an interpolant in that space. Every interpolant in the sequence is constructed using this strategy, which circumvents the Knapsack problem and allows us to restrict our attention to the selection of the optimal polynomial spaces.

The quasi-optimal polynomial space associated with Legendre coefficients is a total degree space with a small logarithmic correction\cite{Beck2014732,TranQuasiOptimal}. However, while the Legendre space is optimal with respect to projection, in the context of interpolation, the quasi-optimal estimate does not account for the effect of the Lebesgue constant. Using estimates of the operator norm of the one dimensional rules, we add a strong correction to the total degree space to arrive at a an estimate for the quasi-optimal interpolation space. Our estimate is parametrized by two vectors associated with the size of the analytic region of the function and the growth of the Lebesgue constant of the interpolation rules.

In order to keep our approach free from a priori assumptions, we present a procedure for dynamically estimating the two vector parameters. For each interpolant in the sequence, we consider the orthogonal decomposition of the interpolant into a linear combination of multivariate Legendre polynomials. Then, we seek the vectors that give the best fit of our quasi-optimal estimate to the decay rate of the Legendre coefficients, i.e., using least-squares approach. The polynomial space used for the construction of the next interpolant in the sequence is optimal with respect to the parameters inferred from the previous interpolant. The number of additional nodes in each interpolant can be chosen arbitrarily, however, few nodes result in more frequent update of the parameter vectors which leads to better accuracy, while larger number of nodes allows for greater parallelization.

The procedure for estimating the quasi-optimal polynomial space can be coupled with any approximation strategy that satisfies a mild assumption regarding the growth of the Lebesgue constant. One potential alternative is to use interpolation based on Fekete points, however, even in moderate dimensions, finding those points involves an ill-conditioned and prohibitively expensive problem. Other popular alternatives are the optimization based methods that construct the an approximation based on minimization of $\ell^2$ (e.g., least-squares\cite{doostan2009least}) or $\ell^1$ (e.g., compressed sensing\cite{doostan2011non}) norms. Those methods can be applied to sets of random samples, however, the number of samples needed to construct the approximation always exceeds the cardinality of the optimal polynomial space. We assume that we can choose the abscissas for each samples and we want to exploit the fact that the range of an interpolation operator has exactly the same  degrees of freedom as the number of interpolation nodes. Thus, the sparse grids interpolants are best suited for our context.

The rest of the paper is organized as follows, in \S\ref{sec:qopt} we derive an estimate of the quasi-optimal interpolation space and we present an iterative procedure for generating a sequence of quasi-optimal polynomial spaces. In \S\ref{sec:SG}, we present a strategy for constructing sparse grids operators with minimal number of nodes and we present several one dimensional interpolation rules. In \S\ref{sec:numerical} we present several numerical examples.

\section{Quasi-optimal polynomial space}\label{sec:qopt}

We consider the problem of approximating a multivariate function $f(\bm y) : \Gamma \to \mathbb{R}$, where $\Gamma \subset \mathbb{R}^d$ is a $d$-dimensional hypercube, i.e., $\Gamma = \bigotimes_{k=1}^d \Gamma_k$ and without loss of generality we let $\Gamma_k = [-1,1]$. We assume that $f(\bm y)$ admits holomorphic extension to a poly-ellipse in complex plane, i.e.,
\begin{assume}[Holomorphic extension]\label{ass:holomorphic}
For a vector $\bm \rho \in \mathbb{R}^d$ with $\rho_k > 1$, the map $\bm z \to f(\bm z)$ is holomorphic in an open neighborhood of the poly-ellipse
\begin{equation}
\mathcal{E}_{\bm \rho} = \bigcup_{\theta \in [0,2\pi]} \mathbin{\mathop{\bigotimes}\displaylimits_{1\leq k \leq d}} \left\{ z_k \in \mathbb{C} : 
		|\Re(z_k)| \leq \frac{\rho_k + \rho_k^{-1}}{2} \cos(\theta), 
		|\Im(z_k)| \leq \frac{\rho_k - \rho_k^{-1}}{2} \sin(\theta)
	 \right\}
\end{equation}
where $\Re(z_k)$ and $\Im(z_k)$ indicate the real and complex part of $z_k$.
\end{assume}
Due to this assumption, we aim at approximaiting $f(\bm y)$ with globally defined polynomials over $\Gamma$. To achieve this goal, we introduce a multivariate polynomial space
\begin{equation}
\nonumber
	\mathcal{P}_{\Lambda(p)}(\Gamma) = span \left\{ \bm y^{\bm \nu} : \bm \nu \in \Lambda(p) \right\},
\end{equation}
in which it will be convenient to use multi-index notation\footnote{For the remainder of the paper we let $\mathbb{N}$ be the set of natural numbers including zero, and $\Lambda, \Theta \subset N^d$ will denote set of multi-indexes. For any two vectors, we define $\bm y^{\bm \nu} = \prod_{k=1}^d y_k^{\nu_k}$ with the usual convention $0^0 = 1$.}, where $\Lambda(p)$ is a sequence of lower multi-index set\footnote{A set $\Lambda$ is caller lower or admissible if $\bm \nu \in \Lambda$ implies $\{ \bm i \in N^d : \bm i \leq \bm \nu \} \subset \Lambda$, where $\bm i \leq \bm \nu$ if and only if $i_k \leq \nu_k$ for all $1 \leq k \leq d$.}. A global polynomial approximation of $f(\bm y)$ in $\mathcal{P}_{\Lambda(p)}(\Gamma)$ has the form
\begin{equation}
\nonumber
	f_{\Lambda(p)} = \sum_{\bm \nu \in \Lambda(p)} c_{\bm \nu} \phi_{\bm \nu}(\bm y),
\end{equation}
where $span\{ \phi_{\bm \nu}(\bm y) : \bm \nu \in \Lambda(p) \} = \mathcal{P}_{\Lambda(p)}(\Gamma)$, and the choice of $\phi_{\bm \nu}(\bm y)$ and $c_{\bm \nu}$ is method specific. To alleviate the curse of dimensionality, the polynomial space $\mathcal{P}_{\Lambda(p)}(\Gamma)$ should be chosen with as few degrees of freedom necessary to approximate $f(\bm y)$ with sufficient accuracy. Using Assumption \ref{ass:holomorphic}, let $\alpha_k = \log( \rho_k )$ for $1 \leq k \leq d$ dictate the anisotropy in each direction, then the most common choice of the polynomial space is the tensor product space $\mathcal{P}_{\Lambda^{TP}(p)}(\Gamma)$, where
\begin{equation}
\nonumber
	\Lambda^{TP}(p) = \{ \bm \nu \in \mathbb{N}^d : \max_{ 1 \leq k \leq d} \alpha_k \nu_k \leq p \},
\end{equation}
and the cardinality depends exponentially on the dimension. Several alternative polynomial spaces have been proposed, namely total degree space with $\Lambda^{TD}(p) = \{ \bm \nu \in \mathbb{N}^d : \bm \alpha \cdot \bm \nu \leq p \}$, hyperbolic cross section $\Lambda^{HC}(p) = \{ \bm \nu \in \mathbb{N}^d : ( \bm \nu + \bm 1 )^{\bm \alpha} \leq p \} $ and Smolyak $\Lambda^{Sm}(p) = \{ \bm \nu \in \mathbb{N}^d : \bm \alpha \cdot \log_2( \bm \nu + \bm 1 ) \leq p \}$, where $\cdot$ indicates vector dot product and $\log( \bm \nu + 1 ) = \bigotimes_{k=1}^d \log( \nu_k + 1 )$ (see \cite{gunzburger2014stochastic} and references therein).

Our approach is motivated by recent work on best $M$-term quasi-optimal Galerkin approximation\cite{Beck2014732,cohen2010convergence,todor2007convergence}. Consider the orthogonal decomposition of $f(\bm y)$
\begin{equation}
\nonumber
	f(\bm y) = \sum_{\bm \nu \in \mathbb{N}^d} c_{\bm \nu} L_{\bm \nu}(\bm y),
\end{equation}
where $L_{\bm \nu}(\bm y)$ are the multivariate Legendre polynomials. For any lower set $\Lambda(p)$ the projection of $f(\bm y)$ onto $\mathcal{P}_{\Lambda(p)}(\Gamma)$ is given by $f_{\Lambda(p)}(\bm y) = \sum_{\bm \nu \in \Lambda(p)} c_{\bm \nu} L_{\bm \nu}(\bm y)$ and the approximation error is
\begin{equation}
\nonumber
	\| f - f_{\Lambda(p)} \|_{L^2}^2 = \sum_{ \bm \nu \not\in \Lambda(p) } | c_{\bm \nu} |^2
\end{equation}
Therefore, the best $M$-term space for projection is the space associated with the $M$ largest coefficients. The coefficients are not known a priori, however, when Assumption \ref{ass:holomorphic} is satisfied, a sharp upper bound to $|c_{\bm \nu}|$ is given by
\begin{equation}
\label{opt:qopt}
	| c_{\bm \nu} | \leq C \cdot exp( - \bm \alpha \cdot \bm \nu ) \prod_{k=1}^d \sqrt{ 2 \nu_k + 1 },
\end{equation}
for some constant $C$ \cite{TranQuasiOptimal,Beck2014732}. The quasi-optimal projection space is then associated with the multi-indexes that yield the largest upper estimate, i.e., 
\begin{equation}
\nonumber
	\Lambda^{\bm \alpha}(p) = \left\{ \bm \nu \in \mathbb{N}^d : \bm \alpha \cdot \bm \nu - \frac{1}{2} \sum_{k=1}^d \log( \nu_k + 0.5 ) \leq p \right\}
\end{equation}
where we derive the condition by taking the log of the right hand side of (\ref{opt:qopt}), changing the sign and ignoring the constant. In this context, $p \in \mathbb{N}$ is an arbitrary variable used to discretize the index space into levels.

\subsection{Quasi-optimal interpolation}

We are interested in constructing approximation using a computationally cheap sampling scheme. Projection results in optimal $L^2(\Gamma)$ error, however, interpolatory approximations are not optimal. Let $\Lambda(L)$ be a lower set and let $f_{\Lambda(L)}(\bm y)$ be in interpolatory approximation of $f(\bm y)$ in $\mathcal{P}_{\Lambda(L)}(\Gamma)$, then
\begin{equation}
\label{opt:interpError}
	\| f - f_{\Lambda(L)} \|_{L^\infty} \leq ( 1 + \mathcal{C}_{\Lambda(L)} ) \inf_{p \in \mathcal{P}_{\Lambda(L)}(\Gamma)} \| f - p \|_{L^\infty},
\end{equation}
where $\mathcal{C}_{\Lambda(L)}$ is the Lebesgue constant, i.e., norm of the interpolation operator, that manifests as an additional penalty term. Observe that
\begin{equation}
\label{opt:infEstimate}
	\inf_{p \in \mathcal{P}_{\Lambda(L)}(\Gamma)} \| f - p \|_{L^\infty} \leq 
	\left\| f - \sum_{\bm \nu \in \Lambda(L)} c_{\bm \nu}  L_{\bm \nu} \right\|_{L^\infty},
\end{equation}
therefore, the dominant polynomial space associated with the infimum term in (\ref{opt:interpError}) can be (heuristically) approximated by estimate (\ref{opt:qopt}), i.e., the difference comes only from using ${L}^2(\Gamma)$ as opposed to the $L^\infty(\Gamma)$ norm. However, to define a proper quasi-optimal interpolation space, (\ref{opt:qopt}) has to be combined with an estimate of the Lebesgue constant. Here we make the following assumption\footnote{In \S\ref{sec:SG} we derive an estimate for the Lebesgue constant associated with our sparse grids construction and we demonstrate that Assumption \ref{ass:Lebesgue} is indeed satisfied.}:

\begin{assume}[Lebesgue constant]\label{ass:Lebesgue}
Let $\mathcal{C}_{\bm \nu}$ indicate the Lebesgue constant associated with the smallest lower set that contains $\bm \nu$, i.e., $\mathcal{C}_{\bm \nu} = \mathcal{C}_{ \{ \bm j \in \mathbb{N}^d : \bm j \leq \bm \nu \} }$, then we assume that
\begin{equation}
\label{opt:LebesgueConstant}
	\mathcal{C}_{\bm \nu} \leq C_{\gamma} \prod_{k=1}^d ( \nu_k + 1 )^{\gamma_k},
\end{equation}
for some constants $C_{\gamma}$ and $( \gamma_k )_{1 \leq k \leq d}$.
\end{assume}

Multiplying (\ref{opt:qopt}) by (\ref{opt:LebesgueConstant}), and using that $\nu_k + 0.5 < \nu_k + 1$ we arrive at the estimate
\begin{equation}
\label{opt:interpEstimate}
	 \mathcal{C}_{\bm \nu} \cdot C \cdot exp( - \bm \alpha \cdot \bm \nu ) \prod_{k=1}^d \sqrt{ 2 \nu_k + 1 } \leq
	 C \cdot C_{\gamma} \cdot \sqrt{2} \cdot exp( - \bm \alpha \cdot \bm \nu ) \prod_{k=1}^d ( \nu_k + 1 )^{\gamma_k + \frac{1}{2}}
\end{equation}
As before, taking the log of (\ref{opt:interpEstimate}), reversing the sign and ignoring the constants, we define the quasi-optimal interpolation space
\begin{equation}
\label{opt:optInterp}
	\Lambda^{\bm \alpha, \bm \beta}(L) = \left\{ \bm \nu \in \mathbb{N}^d : \bm \alpha \cdot \bm \nu + \bm \beta \cdot \log( \bm \nu + \bm 1 ) \leq L \right\},
\end{equation}
where $\beta_k = - \gamma_k - \frac{1}{2}$. The vectors $\bm \alpha$ and $\bm \beta$ give rise to a discretization of the multi-index space, however, the notion of levels does not transcend the specific $\bm \alpha$ and $\bm \beta$, since (\ref{opt:optInterp}) depends on the scaling of the vector components. Given parameter vectors and a sequence of levels $\{ L_n \}_{n=0}^\infty \in \mathbb{N}^d$, we can construct the corresponding quasi-optimal approximations $f_{\Lambda^{\bm \alpha, \bm \beta}(L_n)}(\bm y)$, however, finding suitable $\bm \alpha$ and $\bm \beta$ is of consideration.

\begin{rem}[Preserving the property lower]\label{rem:preserveLower}
The entries of $\bm \beta$ in (\ref{opt:optInterp}) could be negative, hence, $\Lambda^{\bm \alpha, \bm \beta}(L)$ is not necessarily a lower set. Lower sets have the advantage that the polynomial space associated with the multi-indexes is independent from the hierarchical basis used, e.g., Legendre basis, monomials, Newton polynomials. We are interested in constructing a sequence of lower sets, thus, if $\Lambda^{\bm \alpha, \bm \beta}(L)$ is not a lower set, we replace (\ref{opt:optInterp}) with
\begin{equation}
\nonumber
	\Lambda^{\bm \alpha, \bm \beta}(L) = \bigcup_{ \bm \nu \in \mathbb{N}^d, \bm \alpha \cdot \bm \nu + \bm \beta \cdot \log( \bm \nu + \bm 1 ) \leq L } \{ \bm j \in \mathbb{N}^d : \bm j \leq \bm \nu \}
\end{equation}
\end{rem}

\subsection{Estimating the parameters}

Sharp a priori bounds of the region of analytic extension of $f(\bm y)$ and the corresponding $\bm \alpha$ are
seldom available, likewise estimates of the Lebesgue constant are either overly conservative or the constant fluctuates in a wide range and predicting the effective $\gamma_k$ is very difficult. Thus, we need a procedure to dynamically estimate the \emph{effective} parameters $\bm \alpha$ and $\bm \beta$ that give the best fit of (\ref{opt:interpEstimate}) to the behavior of $f(\bm y)$.

Assume that we have constructed an interpolant $f_{\Lambda(L)}(\bm y)$ for some lower set $\Lambda(L)$. Here $\Lambda(L)$ could be chosen according to (\ref{opt:optInterp}) for some $\bm \alpha$ and $\bm \beta$, or according to total degree, hyperbolic cross section or Smolyak formulas\cite{gunzburger2014stochastic}. Since $f_{\Lambda(L)}(\bm y) \in \mathcal{P}_{\Lambda(L)}$ there are coefficients $\hat c_{\bm \nu}$ for $\bm \nu \in \Lambda(L)$ such that
\begin{equation}
\nonumber
	f_{\Lambda(L)}(\bm y) = \sum_{\bm \nu \in \Lambda(L)} \hat c_{\bm \nu} L_{\bm \nu}(\bm y),
\end{equation}
where $ L_{\bm \nu}(\bm y)$ are the multidimensional Legendre polynomials. By orthogonality, each of the coefficients is
\begin{equation}
\label{opt:hatc}
	\hat c_{\bm \nu} = \int_{ \Gamma} f_{\Lambda(L)}(\bm y)  L_{\bm \nu}(\bm y) d{\bm y},
\end{equation}
where the integral can be computed with a multidimensional quadrature rule and note this can be done without additional evaluations of $f(\bm y)$. Since $L_{\bm \nu}(\bm y)$ and $f_{\Lambda(L)}(\bm y)$ are polynomials, it is sufficient to use a quadrature that can integrate exactly all polynomials in $\mathcal{P}_{2 \Lambda(L)}$.

We assume that $| \hat c_{\bm \nu} |$ decay at a rate guided by (\ref{opt:interpEstimate}) for some $\bm \alpha$ and $\bm \beta$, i.e.,
\begin{equation}
\label{opt:decayModel}
	| \hat c_{\bm \nu} | \propto exp( - \bm \alpha \cdot \bm \nu ) ( \bm \nu + \bm 1 )^{-\bm \beta}
	\quad \Longrightarrow  \quad
	\log( | \hat c_{\bm \nu} |)  \approx -\hat C - \bm \alpha \cdot \bm \nu - \bm \beta \cdot \log( \bm \nu + \bm 1 ),
\end{equation}
for some constant $\hat C$. In (\ref{opt:decayModel}), all $\hat c_{\bm \nu}$ are known and we can solve for $\bm \alpha$, $\bm \beta$ and $\hat C$, however, the decay rate of the coefficients is not monotone and hence we look for the parameters that give the ``best fit''. Here  best is used in $\ell^2$ sense, i.e., we infer approximate $\hat{\bm \alpha}$ and $\hat{\bm \beta}$ from the solution to the convex minimization problem
\begin{equation}
\label{opt:leastSQR}
	\min_{ \bm \alpha, \bm \beta, \hat C } \frac{1}{2} \sum_{\bm \nu \in \Lambda(L)} \Big( \hat C + \bm \alpha \cdot \bm \nu + \bm \beta \cdot \log( \bm \nu + \bm 1 ) + \log( | \hat c_{\bm \nu} | ) \Big)^2
\end{equation}
For sufficiently large $\Lambda(L)$, (\ref{opt:leastSQR}) admits a unique solution. However, the accuracy of the estimated $\hat{\bm \alpha}$ and $\hat{\bm \beta}$ strongly depends on the size of $ \Lambda(L)$, in fact, estimates obtained using this least squares approach are valid only for sets close to $\Lambda(L)$. Furthermore, the constant $\hat C$ is not used in (\ref{opt:optInterp}) and since $\hat C exp( - \bm \alpha \cdot \bm \nu ) ( \bm \nu + \bm 1 )^{-\bm \beta}$ does not give an upper bound on the coefficients, $\hat C$ cannot be used to estimate the true approximation error; in this context, $\hat C$ plays the role of a \emph{dummy} variable.

\begin{rem}[Ad hoc stability constraint]\label{rem:AdHoc}
In general, $| \hat c_{\bm \nu} |$ do not decay monotonically and if $\Lambda(L)$ is small, some of the estimated parameters in $\hat{\bm \alpha}$ could be negative. The entries of $\hat{\bm \alpha}$ are associated with the convergence rate in different direction and negative entries indicate that for this specific direction and specific $ \Lambda(L)$ we are observing divergence. In case of a negative $\hat \alpha_k$, we use an \emph{ad hoc} correction, where we replace the negative entries with the smallest positive one, i.e., we assume that in the limit the diverging direction will in fact converge, albeit slowly. 
\end{rem}

The least-squares approach allows us to construct a sequence of polynomial spaces indexed by sets $\{ \Lambda_n \}_{n=0}^\infty$. We start with either $\Lambda_0 = \Lambda^{\bm \alpha_0, \bm \beta_0}(L_0)$ for some initial guess of $\bm \alpha_0$, $\bm \beta_0$ and $L_0$, or $\Lambda_0$ can be chosen according to total degree, hyperbolic cross section or Smolyak formulas. Then from $\Lambda_n$ and the interpolant $f_{\Lambda_n}(\bm y)$, we infer the best fit parameters $ \hat{\bm \alpha}_{n+1}$ and $\hat{\bm \beta}_{n+1}$ and construct the next index set as
\begin{equation}
\label{opt:LambdaSequence}
	\Lambda_{n+1} = \Lambda_n \bigcup \Lambda^{ \hat{\bm \alpha}_{n+1}, \hat{\bm \beta}_{n+1} }( L_{n+1} ),
\end{equation}
where we take the union to ensure that the sequence is nested and $L_{n+1}$ controls the number of additional degrees of freedom in the polynomial space. Using small $L_{n+1}$, i.e., taking the smallest $L_{n+1}$ such that $\Lambda^{ \hat{\bm \alpha}_{n+1}, \hat{\bm \beta}_{n+1} }( L_{n+1} ) \not\subset \Lambda_n$, leads to a more frequent update of the parameters $\hat{\bm \alpha}$ and $\hat{\bm \beta}$ and hence better accuracy, however, larger $L_{n+1}$ leads to more samples needed for constructing $f_{\Lambda_{n+1}}(\bm y)$ and hence more opportunity for parallelization.

Next, we present a specific sparse grids strategy for constructing $f_{\Lambda_n}(\bm y)$.

\section{Sparse grids interpolation}\label{sec:SG}

In this section we present a general sparse grids interpolation approach, which consists of evaluating $f(\bm y)$ at a set of nodes $\bm y_1, \cdots, \bm y_m \in \Gamma$ and constructing an interpolant $f_{\Lambda(L)}(\bm y) \in \mathcal{P}_{\Lambda(L)}$, where $\Lambda(L)$ is an arbitrary lower index set, i.e., not necessarily constructed according to (\ref{opt:optInterp}). The properties of the SG interpolant are determined by the Lebesgue constant and growth of nodes in the one dimensional family of interpolants that induce the grid. Minimizing the number of nodes is desirable and ideally we want the number of samples to not exceed the cardinality of $\Lambda(L)$, however, interpolants with more nodes often times result in smaller operator norm and potentially more accurate approximation. For a given $\Lambda(L)$ and one dimensional rule, we present a strategy for constructing the SG with smallest number of nodes that produces an interpolant in $\mathcal{P}_{\Lambda(L)}$. We also present several novel one dimensional rules.

\subsection{Constructing optimal interpolant} A nested one dimensional family of interpolation rules is defined by a distinct sequence of nodes $\{ y_j \}_{j=1}^{\infty} \in [-1,1]$ and a strictly increasing growth function $m: \mathbb{N} \to \mathbb{N}$, i.e., $m(0) > 0$ and $m(l+1) > m(l)$. For $l = 0, 1, 2, \cdots$ we define the $l$-th level interpolation operator
\begin{equation}
\nonumber
	\mathcal{U}^{m(l)} : \mathbb{C}^0(\Gamma) \to \mathcal{P}_{m(l)-1}([-1,1]), \qquad \mbox{by} \qquad
	\mathcal{U}^{m_l}[g](y) = \sum_{j=1}^{m(l)} g(y_j) \psi^l_j(y),
\end{equation}
where $g(y) \in \mathbb{C}^0([-1,1])$, the Lagrange basis functions are $\psi_j^l(y) = \prod_{ i = 1, i \neq j }^{m_l} \frac{y - y_{i}}{y_j - y_{i}}$ and $\mathcal{P}_{m(l)-1}([-1,1]) = span\{ y^\nu : 0 \leq \nu \leq m(l) - 1\}$. In addition, for $l>0$, we define the surplus operators $ \Delta^{m(l)} = \mathcal{U}^{m(l)} - \mathcal{U}^{m(l-1)}$ and for notational convenience let $\Delta^0 = \mathcal{U}^{m(0)}$. Several specific examples of $y_j$ and $m_l$ are listed in Table \ref{table:rules} in \S\ref{sec:1Drules}. Note, we are explicitly assuming the interpolants have nested nodes and strictly increasing $m(l)$, see Remarks \ref{rem:Nested} and \ref{rem:growth}.

Taking the $d$-dimensional tensors, we have
\begin{equation}
\nonumber
	\bm m(\bm i) = \bigotimes_{k=1}^d m(i_k) : \mathbb{N}^d \to \mathbb{N}^d, \qquad
	\bm y_{\bm j} = \bigotimes_{k=1}^d y_{j_k} \in \Gamma, \qquad
	\Psi_{\bm j}^{\bm i}(\bm y) = \prod_{k=1}^d \psi_{j_k}^{i_k},
\end{equation}
where $\bm i, \bm j \in \mathbb{N}^d$, $\bm y \in \Gamma$ and $ \psi_{j_k}^{i_k} $ is evaluated at the corresponding $k$-th component of $\bm y$. The tenor operators are given by
\begin{equation}
\nonumber
	\Delta^{\bm m(\bm i)} = \bigotimes_{k=1}^d \Delta^{m(i_k)}, \qquad
	\mathcal{U}^{\bm m(\bm i)} = \bigotimes_{k=1}^d \mathcal{U}^{m(i_k)}, \qquad %= \sum_{ \bm j \leq \bm i } \Delta^{\bm m(\bm j)}
	\mathcal{U}^{\bm m(\bm i)}[f](\bm y) = \sum_{ \bm 1 \leq \bm j \leq \bm m(\bm i)} f(\bm y_{\bm j}) \Psi_{\bm j}^{\bm i}(\bm y)
\end{equation}
Note that by the telescoping properties of $\Delta^{\bm m(\bm i)}$ we have that $\mathcal{U}^{\bm m(\bm i)} = \sum_{ \bm j \leq \bm i } \Delta^{\bm m(\bm j)}$, i.e., each full tensor operator $\mathcal{U}^{\bm m(\bm i)}$ can be decomposed into a sum of surplus operators.

%%%%%%%%%%%%%%%%%%%%%%%%%%%%%%%%%%%%%%%%%%%%%%%%%%%%%%%%%%%%%%%%%%%%%
%
%    PARAGRAPH : Delta surplus operator
%
%%%%%%%%%%%%%%%%%%%%%%%%%%%%%%%%%%%%%%%%%%%%%%%%%%%%%%%%%%%%%%%%%%%%%

Every multidimensional polynomial space can be included in the range of a full tensor interpolation operator, however, full tensors have a rigid structure and often times require an excessive number of samples. Sparse grids offer a flexible alternative, where the interpolant is constructed from a sparse set of the surplus operators $\T(L)$, i.e.,
\begin{equation}
\label{sg:mDGenInterp}
	I^m_{\T(L)} = \sum_{\bm i \in \T(L)} \Delta^{\bm m(\bm i)}
\end{equation}
For any lower index set $\T(L)$, (\ref{sg:mDGenInterp}) is an interpolation operator with nodes $\{ \bm y_{\bm j} \}_{\bm j \in \T_m(L)}$, where
\begin{equation}
\label{sg:nodes}
	\T_m(L) = \bigcup_{\bm i \in  \T(L)} \{ \bm j \in \mathbb{N}^d : \bm 1 \leq \bm j \leq \bm m(\bm i) \}
		 = \bigcup_{\bm i \in  \T(L)} \{ \bm j \in \mathbb{N}^d : \bm m_{\bm i - \bm 1 } + \bm 1 \leq \bm j \leq \bm m(\bm i) \}
\end{equation}
and $I^m_{\T(L)}[f] \in \mathcal{P}_{\T^m(L) - \bm 1}$ with $ \T_m(L) - \bm 1 = \{ \bm j \in \mathbb{N}^d : \bm j + \bm 1 \in \T_m(L) \}$ \cite{NobileQuasiOptimalInterpolation}. By definition of $\Delta^{\bm m(\bm i)}$, there exists a set of integer weights $\{t_{\bm j} \}_{\bm j \in \T(L)}$, satisfying the system of equations $\sum_{\bm j \leq \bm i, \bm j \in \T(L)} t_{\bm j} = 1$ for every $\bm i \in \T(L)$, and the interpolant can be written explicitly as
\begin{equation}
\label{sg:explicitI}
	I^m_{\T(L)}[f](\bm y) = \sum_{ \bm j \in \T_m(L) } f(\bm y_{\bm j}) \sum_{ \bm i \in \T(L), \bm m(\bm i) \leq \bm j } t_{\bm i} \Psi^{\bm i}_{\bm j} (\bm y)
\end{equation}
Interpolant (\ref{sg:explicitI}) is constructed from the samples $f(\bm y_{\bm j})$ for $\bm j \in \T_m(L)$, and since the sets of the second union of (\ref{sg:nodes}) are disjoint, we have a direct relationship between $\T(L)$, $m(l)$ and the number of nodes.

\begin{rem}\label{rem:Nested} If the one dimensional family of rules is not nested then generally (\ref{sg:mDGenInterp}) is not an interpolant. Even in the non-nested case, operator of the form (\ref{sg:mDGenInterp}) can produce an accurate approximation to $f(\bm y)$\cite{nobile2008sparse,NobileQuasiOptimalInterpolation}, however, the approximation belongs to a space of polynomials with cardinality much less than the required number of function evaluations. Excess sampling unnecessarily increases the computational cost and thus we restrict our attention to nested rules.
\end{rem}

Next, we present the main result of this section, where we consider the optimal choice of tensor set that will guarantee the range of the interpolation operator includes a given polynomial space.

%%%%%%%%%%%%%%%%%%%%%%%%%%%%%%%%%%%%%%%%%%%%%%%%%%%%%%%%%%%%%%%%%%%%%
%
%    THEOREM : Minimal interpolation
%
%%%%%%%%%%%%%%%%%%%%%%%%%%%%%%%%%%%%%%%%%%%%%%%%%%%%%%%%%%%%%%%%%%%%%

\begin{theorem}[Minimal polynomial interpolant]\label{theorem:opt}
Let $ \Lambda(L)$ be a lower set of polynomial multi-indexes and define
\begin{equation}
\label{sg:topt}
	 \T^{opt}(L) = \left\{ \bm i \in \mathbb{N}^d : \bm m(\bm i - \bm 1) \in  \Lambda(L) \right\},
\end{equation}
where for notational convenience we let $m(-1) = 0$. Then, $\mathcal{P}_{ \Lambda(L)} \subset \mathcal{P}_{\T^{opt}_m(L) - \bm 1}$ and with respect to the particular $m(l)$, $ \T^{opt}(L)$ is minimal, i.e., if $ \T(L)$ is a lower set such that $\mathcal{P}_{ \Lambda(L)} \subset \mathcal{P}_{\T_m(L) - \bm 1}$, then $ \T(L)$ is a superset of $ \T^{opt}(L) $.
\end{theorem}
\emph{Proof:} For arbitrary $\bm j \in  \Lambda(L)$, by monotonicity of $m(l)$, there is $\bm i$ such that $m(i_k-1) \leq j_k \leq m(i_k)-1$ for $k=1,2,\cdots,d$. Then, according to (\ref{sg:topt})
\begin{equation}
\nonumber
	\bm m(\bm i - \bm 1) \leq \bm j \Longrightarrow \bm i \in  \T^{opt}(L), \quad \mbox{and} \quad
	\bm j \leq \bm m(\bm i) - \bm 1 \Longrightarrow \bm j \in  \T^{opt}_m(L) - \bm 1.
\end{equation}

To show that $ \T^{opt}(L) $ is minimal, suppose $ \T(L)$ is a lower set and $ \Lambda(L) \subset \T_m(L) - \bm 1$. For arbitrary $\bm i \in \T^{opt}(L)$, let $\bm j = \bm m(\bm i -\bm 1)$, then
\begin{equation}
\nonumber
	\bm j = \bm m(\bm i -\bm 1) \Longrightarrow \bm j \in  \Lambda(L) \Longrightarrow \bm j \in  \T_m(L) - \bm 1
	\Longrightarrow \mbox{exists } \bm i' \in  \T(L) \mbox{ s.t. } \bm j \leq \bm m(\bm i') - \bm 1.
\end{equation}
Thus, $\bm m(\bm i - \bm 1) \leq \bm m(\bm i') - \bm 1$ and by monotonicity of $m(l)$ we have $\bm i \leq \bm i'$. Since $\bm i' \in  \T(L)$ and $ \T(L)$ is a lower set, $\bm i \in \T(L)$.  \hfill\ensuremath{\square} %\mbox{\square}

\begin{cor}\label{cor:opt} For $ \Lambda^{\bm \alpha, \bm \beta}(L)$ defined in (\ref{opt:optInterp}), the optimal tensor set is
\begin{equation}
\label{sg:talphabeta}
	\T^{\bm \alpha, \bm \beta}(L) = \left\{ \bm i \in \mathbb{N}^d : \bm \alpha \cdot \bm m(\bm i - \bm 1) + \bm \beta \cdot \log( \bm m(\bm i - \bm 1) + \bm 1 ) \leq L \right\}
\end{equation}
\end{cor}

Armed with this result, we define the sparse grids inteprolants associated with the sequence of polynomial spaces defined in (\ref{opt:LambdaSequence}). Let $\T_0$ be the optimal tensor set associated with $\Lambda_0$, then $f_{\Lambda_0}(\bm y) = I_{\T_0}^m[f](\bm y)$ and for $n \in \mathbb{N}$
\begin{equation}
\nonumber
	\T_{n+1} = \T_{n} \bigcup \T^{\hat{\bm \alpha}_{n+1}, \hat{\bm \beta}_{n+1}}(L_{n+1}), 
	\qquad \mbox{and} \qquad
	f_{\Lambda_{n+1}}(\bm y) = I_{\T_{n+1}}[f](\bm y)
\end{equation}
where $\hat{\bm \alpha}_{n+1}$ and $\hat{\bm \beta}_{n+1}$ are estimated from (\ref{opt:leastSQR}). Note that unless $m(l) = l+1$ the interpolants constructed in this fashion will be associated with polynomial spaces larger than $\mathcal{P}_{\Lambda_n}$. The traditional Knapsack approach uses $m(l)$ explicitly into the definition of the profit index\cite{NobileQuasiOptimalInterpolation,bungartz2004sparse}, however, in our context, rules with $m(l) > l+1$ are associated with smaller Lebesgue constant and thus $m(l)$ affects $\Lambda_n$ implicitly through the $\bm \beta$ term.

\begin{rem} \label{rem:growth} Isotropic total degree space defined by $\sum_{k=1}^d j_k \leq L$ is an example of a particular polynomial space of interest, in fact, large amount of the initial work in sparse grids was aimed at constructing total degree interpolants. The work\cite{novak1999simple} give an optimal construction using one dimensional rules with occasionally repeating number of nodes (as opposed to strictly increasing). However, Corollary \ref{cor:opt} with $\bm \alpha = \bm 1$ and $\bm \beta = \bm 0$ gives us the same result as the slow-growth method and hence the assumption $m(l) < m(l+1)$ is not restrictive.
\end{rem}

%%%%%%%%%%%%%%%%%%%%%%%%%%%%%%%%%%%%%%%%%%%%%%%%%%%%%%%%%%%%%%%%%%%%%
%
%    SUBSECTION : Lebesgue ``constant''
%
%%%%%%%%%%%%%%%%%%%%%%%%%%%%%%%%%%%%%%%%%%%%%%%%%%%%%%%%%%%%%%%%%%%%%

\subsection{Lebesgue constant}\label{subsec:LebesgueConstant} 

In the one dimensional context, the norm of $\mathcal{U}^{m(l)}$ can be estimated numerically and, in some cases, sharp theoretical estimates are also available. Thus, we assume there exists $\{ \lambda_l \}_{l=0}^\infty \in \mathbb{R}$ so that $\| \mathcal{U}^{m(l)} \|_{L^\infty} \leq \lambda_l$, and for specific examples see \S\ref{sec:1Drules}. However, even if $\lambda_l$ is sharp, there is no known \emph{sharp} analytic estimate of $\| I_{ \T(L)}^m \|_{L^\infty}$ for a general $ \T(L)$ and numerical estimates are computationally impractical. In the case when $\lambda_l$ exhibits polynomial growth, i.e., $\lambda_l \leq C_\gamma (l+1)^\gamma$ for some $C_\gamma, \gamma \in \mathbb{R}$, Lemma 3.1 in \cite{chkifa2014high} shows that
\begin{equation}
\label{sg:mDGenNorm}
	\| I_{ \T(L)}^m \|_{\mathbb{C}^0} \leq C^d_\gamma ( \# \T(L) )^{\gamma+1}
\end{equation}
where $ \# \T(L) $ indicates the umber of elements in $\T(L)$. This result is not sharp and the effective operator norm is usually much smaller than what (\ref{sg:mDGenNorm}) suggests, however, estimate (\ref{sg:mDGenNorm}) indicates that the two main factors contributing to the operator norm are the growth of $\lambda_l$ (i.e., $C_\gamma$ and $\gamma$) and the cardinality of $ \T(L)$. 

The polynomial growth of $\lambda_l$ is sufficient to satisfy Assumption \ref{ass:Lebesgue}. The norm of a full tensor operator is $\| \mathcal{U}^{\bm m( \bm i)} \|_{L^\infty} \leq \prod_{k=1}^d \lambda_{i_k}$, and the associated polynomial space is indexed by $\{ \bm \nu \in \mathbb{N}^d : \bm \nu \leq \bm m( \bm i ) - \bm 1 \}$. By monotonicity of $m(l)$, if $\bm i$ is the index of the smallest full tensor containing $\bm \nu$, then $\bm i \leq \bm \nu$ and
\begin{equation}
\label{sg:assump2}
	\mathcal{C}_{\bm \nu} \leq \prod_{k=1}^d \lambda_{i_k} \leq C^d_\gamma \prod_{k=1}^d ( i_k + 1 )^\gamma
	\leq C^d_\gamma \prod_{k=1}^d ( \nu_k + 1 )^\gamma
\end{equation}
Note that (\ref{sg:assump2}) implies that $\gamma_k$ is independent from $k$, however, a one dimensional family of rules can exhibit very slow increase of the Lebesgue constant for the first several nodes, followed by much sharper increase. Thus, when $f(\bm y)$ exhibits strong anisotropic behavior (i.e., the components of $\bm \alpha$ differ significantly), there is a large discrepancy in the number of one dimensional nodes associated with various directions, and different $\gamma_k$  yield a much sharper estimate.

%%%%%%%%%%%%%%%%%%%%%%%%%%%%%%%%%%%%%%%%%%%%%%%%%%%%%%%%%%%%%%%%%%%%%%%%%%%%%%%%%%%%%%%%%%%%%%%%%%%%%%%%%%%%%%%%%%%%%
%
%  SUBSECTION: 1-D Rules
%
%%%%%%%%%%%%%%%%%%%%%%%%%%%%%%%%%%%%%%%%%%%%%%%%%%%%%%%%%%%%%%%%%%%%%%%%%%%%%%%%%%%%%%%%%%%%%%%%%%%%%%%%%%%%%%%%%%%%%

\subsection{One dimensional rules} \label{sec:1Drules}

\begin{figure}
\includegraphics[scale=0.5]{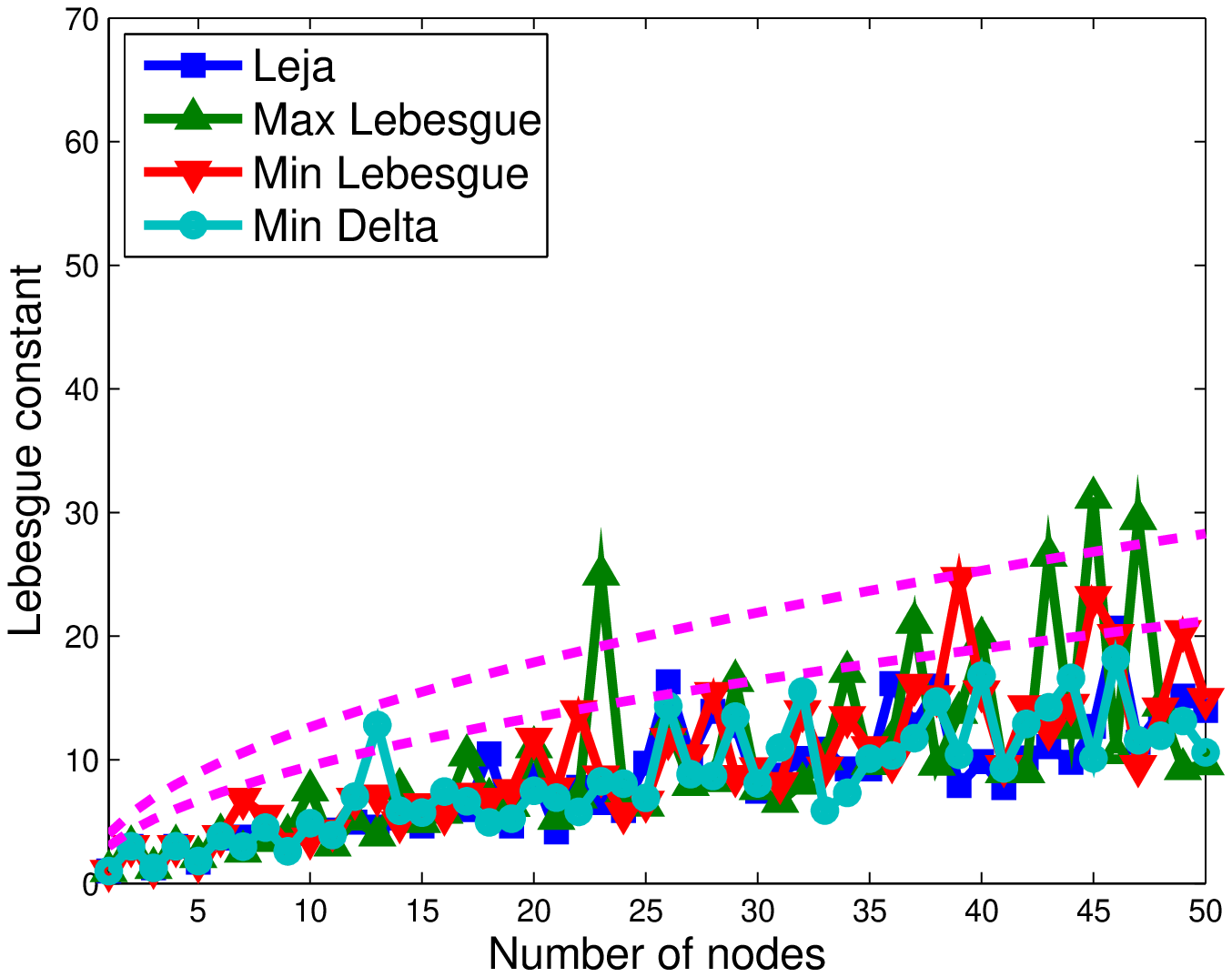}
\includegraphics[scale=0.5]{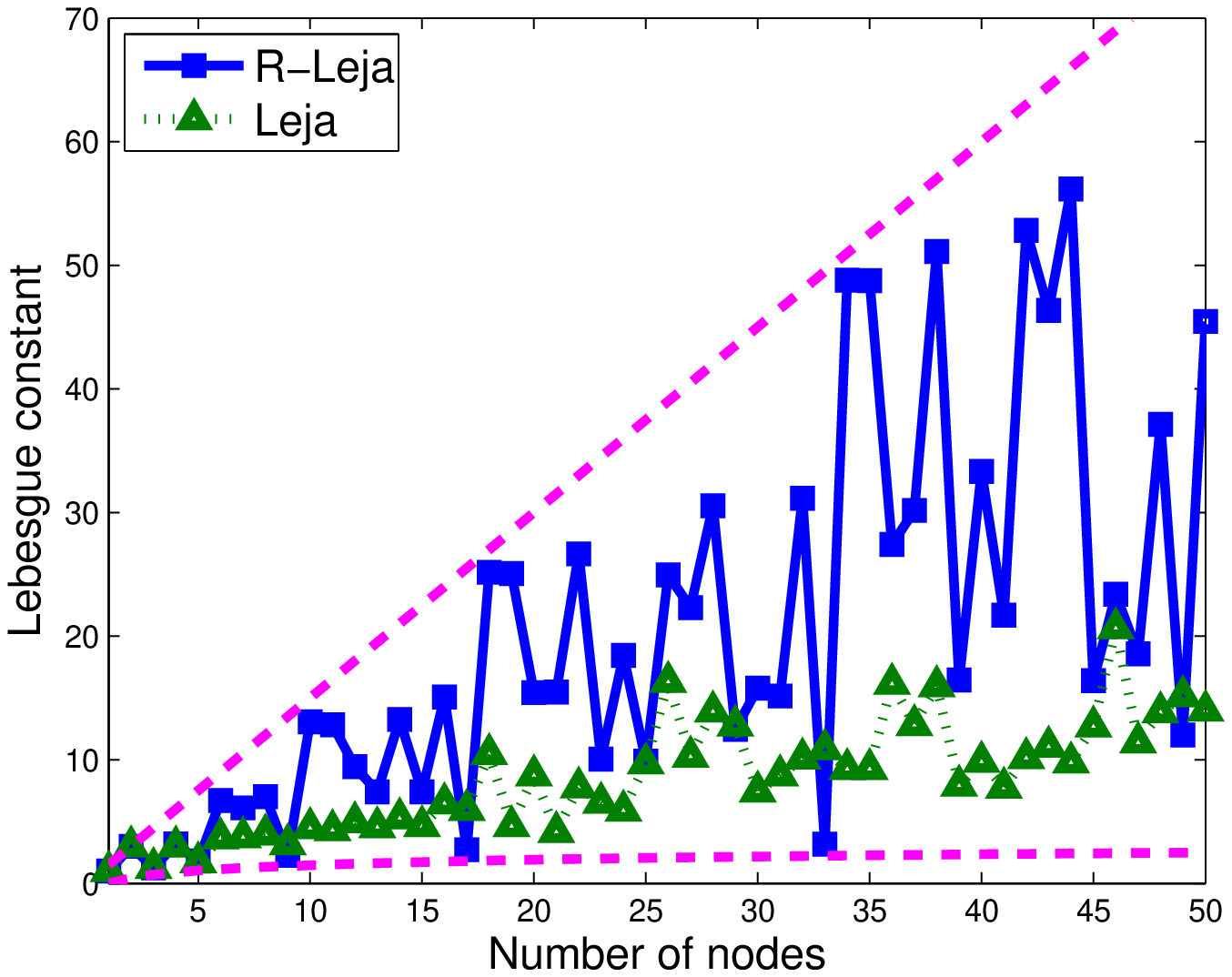}
\caption{Interpolation operator norm for different sequence of points. Left: the two dashed lines correspond to $3\sqrt{l+1}$ and $4\sqrt{l+1}$. Right: the reference lines are $\frac{3}{2} (l+1)$ and $\frac{2}{\pi}\log( 2^l+1 )$. } \label{fig:LebPenalty}
\end{figure}

In this section, we present multiple one dimensional interpolation rules with different nodes $y_j$, growth sequence $m_l$ and operator norm $\lambda_l$. Ideally, we want a rule with slowly growing $m(l)$ and $\lambda_l$, since small $m_l$ results in fewer interpolation nodes and smaller $\lambda_l$ leads to better accuracy of the approximation, however, those are competing goals. Rules with rapidly increasing $m(l)$ lead to $ \T_n$ with small cardinality and interpolants with small Lebesgue constant. However, rapid growth of $m(l)$ also leads to interpolants with more nodes than degrees of freedom of $\Lambda_n$. It is hard to predict a priori the optimal relation between $m(l)$ and $\lambda_l$, and in what follows, we present a list of several candidate interpolation rules.

The roots and extrema of Chebyshev polynomials are a popular choice for interpolation nodes and the Chenshaw-Curtis\cite{clenshaw1960method} rule is one of the most widely used. The Chenshaw-Curtis nodes $y_j$ are defined by
\begin{equation}
\label{eq:CCnodes}
	y_1 = 0, \quad y_2 = 1, \quad y_3 = -1, \qquad \mbox{for } j> 3, \quad y_j = \cos \Big( 2^{ -\lceil \log_2(j-1) \rceil } ( 2 j - 3 ) \pi \Big),
\end{equation}
where $\lceil x \rceil = \min \{ i \in \mathbb{Z} : i \geq x \}$. The growth sequence starts with $m(0) = 1$ and for $l>0$ we have $m(l) = 2^l + 1$. The operator norm increases logarithmically in number of nodes
\begin{equation}
\label{eq:ccLConst}
	\lambda_l = \frac{2}{\pi} \log( m(l) -1 ) + 1 = \frac{2}{\pi} \log( 2^l ) + 1.
\end{equation}
Another similar example is the Fejer type-2 interpolation based on the interior roots of Chebyshev polynomials\cite{fejer1933infinite}. The nodes are given by
\begin{equation}
\label{eq:fejer2nodes}
	y_j = \cos \Big( 2^{ -\lceil \log_2(j+1) \rceil } ( 2 j + 1 ) \pi \Big),
\end{equation}
with growth sequence $m(l) = 2^{l+1}-1$. Similar to Clenshaw-Curtis, the operator norm exhibits logarithmic dependence on $m(l)$.

\begin{rem}%[Logarithmic increase of operator norm]
For interpolation rules that exhibit logarithmic increase of $\lambda_l$, replacing $( \nu_k + 1 )^{\gamma}$ in Assumption \ref{ass:Lebesgue} with $\log( \nu_k + 1 )$ may yield a sharper estimate. However, the effects of the logarithmic term are negligible even for relatively small $\bm \alpha$. See Figure \ref{fig:KL7D_selection} in \S\ref{sec:numerical}.
\end{rem}

A nested rule based on Chebyshev nodes and a slower growing $m_l$ is the $\mathcal{R}$-Leja sequence presented in \cite{chkifa2013lebesgue}. Let $\{ \theta_k \}_{k=1}^{\infty}$ be the sequence defined recursively by
\begin{equation}
\label{eq:rlejaTheta}
	\theta_1 = 0, \quad \theta_2 = \pi, \quad \theta_3 = \frac{\pi}{2}, \quad \mbox{for $j>3$, }
		\theta_j = \left\{ \begin{array}{ll}
			\theta_{j-1} + \pi, & \mbox{$j$ is odd}\\
			\frac{1}{2}\theta_{\frac{j}{2}+1}, & \mbox{$j$ is even}\\
		\end{array} \right.
\end{equation}
then the $\mathcal{R}$-Leja nodes are $y_j = \cos( \theta_j )$ with $m(l) = l+1$, and $\lambda_l \leq 2(l+1)^2 \log(l+1)$. The single point growth gives great flexibility, since the number of nodes needed to construct $f_{\Lambda_n}(\bm y)$ always equals the number of indexes in $\Lambda_n$, however, the operator norm of $\mathcal{R}$-Leja based interpolants is usually larger than Clenshaw-Curtis and Fejer rules.

The quadratic estimate of the Lebesgue constant for the $\mathcal{R}$-Leja sequence is sharp for the worst case, however, the actual penalty fluctuates between quadratic and logarithmic (in the number of nodes), which gives rise to a family of rules with different $m(l)$\cite{chkifa2013lebesgue}. First, we define the centered $\mathcal{R}$-Leja sequence as
\begin{equation}
\label{eq:centeredRleja}
	y_1 = 0, \quad y_2 = 1, \quad y_3 = -1, \quad y_j = \cos( \theta_j ),
\end{equation}
where $\theta_j$ is defined as in (\ref{eq:rlejaTheta}). With the centered rule, if we select $m(l) = 2^l + 1$ (and $m(0) = 1$), then the resulting one dimensional interpolants are identical to those defined by the Clenshaw-Curtis rule. In general, $\mathcal{R}$-Leja sub-sequences with exponentially growing $m(l)$ exhibit linear (in $l$) increase of $\lambda_l$. Two specific examples include the $\mathcal{R}$-Leja double-$2$ rule defined by
\begin{equation}
\label{eq:rlejaDouble2}
	m(0) = 1, \quad m(1) = 3, \quad \mbox{for $l>1$,} \quad m(l) = 2^{ \lfloor \frac{l}{2} \rfloor + 1 } \left( 1 + \frac{l}{2} - \Big\lfloor \frac{l}{2} \Big\rfloor \right) + 1,
\end{equation}
and the $\mathcal{R}$-Leja double-$4$ rule defined by
\begin{equation}
\label{eq:rlejaDouble4}
	m(0) = 1, \quad m(1) = 3, \quad \mbox{for $l>1$,} \quad m(l) = 2^{ 2 + \lfloor \frac{l-2}{4} \rfloor } \left( 1 + \frac{l-2}{4} - \Big\lfloor \frac{l-2}{4} \Big\rfloor \right) + 1,
\end{equation}
where $\lfloor x \rfloor = \max \{ i \in \mathbb{Z} : i \leq x \}$. The word \emph{double} in the name refers to the fact that $l$ appears in the exponents of $2$ (i.e., we are doubling the number of nodes) and the numbers $2$ and $4$ refer to the delay of the doubling (i.e., $m(l+2) - 1 = 2 ( m(l) - 1)$ and $m(l+4) - 1 = 2 ( m(l) - 1 )$ respectively). Yet another option is to use the odd sub-sequence, namely $m(l) = 2l+1$. The $\mathcal{R}$-Leja odd rules result in interpolants with symmetric distribution of nodes and slightly lower operator norm. See Table \ref{table:rules} for a summary of the $\mathcal{R}$-Leja rules.

The $\mathcal{R}$-Leja sequence is constructed as the solution to a greedy optimization problem defined on the unit disk in the complex plane, the resulting complex nodes are projected onto the real line. The optimization problem can also be defined on $[-1,1]$ as $y_1 = 0$ and
\begin{equation}
\label{eq:optleja}
	y_{j+1} = \underset{y \in [-1,1]}{\text{argmax}} \prod_{ i = 1 }^{ j } \big| y - y_i \big|,
\end{equation}
where, if the optimization admits more than one solution, we take the right-most one. This construction leads to Leja interpolation\cite{de2004leja} and the number of points in Leja interpolants can grow one at a time, i.e., $m(l) = l + 1$, or we can use only the odd rules, i.e., $m(l) = 2l+1$. Unlike the $\mathcal{R}$-Leja sequence, there is no known sharp estimate of the operator norm of Leja interpolants, however, numerical tests shows that for $l \leq 50$ we can take $\lambda_l \approx 3 (l+1)^{1/2}$.

Alternative greedy sequences can be constructed by replacing (\ref{eq:optleja}) with maximization of the Lebesgue function
\vspace{-0.3cm}
\begin{equation}
\label{eq:maxLebesgue}
	y_{j+1} = \underset{y \in [-1,1]}{\text{argmax}} \sum_{ j' = 1 }^j \prod_{ i = 1, i \neq j' }^j \Big| \frac{y - y_i}{y_{j'} - y_i} \Big|,
\vspace{-0.3cm}
\end{equation}
minimization of $\| \mathcal{U}^{m(l)} \|_{L^\infty}$
\vspace{-0.4cm}
\begin{equation}
\label{eq:minLebesgue}
	y_{j+1} = \underset{y \in [-1,1]}{\text{argmin}} \max_{ y' \in [-1,1] } \prod_{i=1}^j \Big| \frac{y' - y_i}{y - y_i} \Big| 
		+ \sum_{ j' = 1 }^j \Big| \frac{y' - y}{y_{j'} - y} \Big|  \prod_{ i = 1, i \neq j' }^j \Big| \frac{y' - y_i}{y_{j'} - y_i} \Big|,
\vspace{-0.3cm}
\end{equation}
or minimization of $\| \Delta^{m(l)} \|_{L^\infty}$
\vspace{-0.3cm}
\begin{equation}
\label{eq:minDelta}
	y_{j+1} = \underset{y \in [-1,1]}{\text{argmin}} \max_{ y' \in [-1,1] } 
		\left( 1 + \sum_{i=1}^j \prod_{ j'=1, j' \neq i }^j \Big| \frac{y - y_{j'}}{y_i - y_{j'}} \Big| \right)
		\prod_{j'=1}^j \Big| \frac{y' - y_{j'}}{y - y_{j'}} \Big|.
\vspace{-0.3cm}
\end{equation}
In (\ref{eq:maxLebesgue})-(\ref{eq:minDelta}) the growth can be prescribed either as a one or two points (i.e., $m(l) = l+1$ or $m(l) = 2l+1$). Numerical estimates of $\lambda_l$ for each sequence are given in Figure \ref{fig:LebPenalty}.

\begin{table}
\begin{center}
\begin{tabular}{|l|l|l|l|}
\hline
Name & $m_l$ & $\lambda_l$ & $y_k$ \\
\hline
Clenshaw-Curtis               & $1,2^l+1$   & $\frac{2}{\pi} \log( 2^l + 1)$ & see (\ref{eq:CCnodes}) \\
Fejer type 2                  & $2^{l+1}-1$ & $\frac{2}{\pi} \log( 2^{l+1} - 1 )$ & see (\ref{eq:fejer2nodes}) \\
$\mathcal{R}$-Leja            & $l+1$       & $1.5(l+1)$ & $\cos( \theta_k )$ (see \ref{eq:rlejaTheta}) \\
$\mathcal{R}$-Leja double-$2$ & see (\ref{eq:rlejaDouble2}) & $1.5(l+1)$ & see (\ref{eq:centeredRleja}) \\
$\mathcal{R}$-Leja double-$4$ & see (\ref{eq:rlejaDouble4}) & $1.5(l+1)$ & see (\ref{eq:centeredRleja}) \\
$\mathcal{R}$-Leja odd        & $2l+1$                      & $3(l+1)$ & see (\ref{eq:centeredRleja}) \\
Leja                          & $l+1$  & $3\sqrt{l+1}$ & see (\ref{eq:optleja}) \\
Leja odd                      & $2l+1$ & $6\sqrt{l+1}$ & see (\ref{eq:optleja}) \\
max-Lebesgue                  & $l+1$  & $4\sqrt{l+1}$ & see (\ref{eq:maxLebesgue}) \\
max-Lebesgue odd              & $2l+1$ & $8\sqrt{l+1}$ & see (\ref{eq:maxLebesgue}) \\
min-Lebesgue                  & $l+1$  & $4\sqrt{l+1}$ & see (\ref{eq:minLebesgue}) \\
min-Lebesgue odd              & $2l+1$ & $8\sqrt{l+1}$ & see (\ref{eq:minLebesgue}) \\
min-delta                     & $l+1$  & $3\sqrt{l+1}$ & see (\ref{eq:minDelta}) \\
min-delta odd                 & $2l+1$ & $6\sqrt{l+1}$ & see (\ref{eq:minDelta}) \\
\hline
\end{tabular}
\caption{Summary of one dimensional rules. Note that the increase of the Lebesgue penalty $\lambda_l$ is mesured empirically form the first $50$ nodes of the corresponding sequence.}\label{table:rules}
\end{center}
\end{table}

%\vspace{-1cm}

\begin{rem}[Alternative interpolation rules]
Each of the optimization problems (\ref{eq:optleja}) - (\ref{eq:minDelta}) can be redefined over a multidimensional domain and a greedy procedure can be devised that constructs an interpolant for an arbitrary $ \Lambda(L)$. However, The optimization problem is ill-conditioned and feasible only for few dimensions and moderate cardinality of $\Lambda(L)$, and therefore not practical. A notable exception is the Magic Points procedure\cite{maday2009general}, which is an extension of (\ref{eq:optleja}) to a multidimensional domain of arbitrary geometry, and the method can be used for non-polynomial interpolation. However, in the case of a hypercube $\Gamma$ and a lower polynomial space, the functional associated with the Magic Points greedy problem is a product of one dimensional functionals of type (\ref{eq:optleja}), and the maximum is a tensor of one dimensional Leje nodes. Thus, one possible realization of the Magic Points algorithm is a sparse grid induced from Leja nodes.
\end{rem}

\subsection{Summary of our method} Our dynamically adaptive sparse grids interpolaiton strategy is summarized in the following algorithm:

\begin{algorithm}[Anisotropic Dynamically Adaptive Multidimensional Approximation]\label{alg:ADAMA}
\begin{algorithmic}
\STATE
\STATE Given target function $f(\bm y) \in \mathbb{C}^0(  \Gamma )$, where $ \Gamma$ is a hypercube in $\mathbb{R}^d$
\STATE Select one dimensional nodes $\{ y_j \}_{j=1}^{\infty}$ and growth rule $m(l)$ (e.g., form Table \ref{table:rules})
\STATE Select initial $\Lambda_0$, e.g., (\ref{opt:optInterp}), total degree, hyperbolic cross section or Smolyak formulas.
\STATE Let $\T_0 = \T^{opt}_0$ according to (\ref{sg:topt})
\FOR{ $n = 0,1,2,\cdots$ \TO $\infty$ }
\STATE Construct $f_{\Lambda_n}(\bm y) = I_{\T_n}[f](\bm y)$ according to (\ref{sg:explicitI}), i.e., compute the necessary $f(\bm y_{\bm j})$
\STATE Compute the coefficients $\hat c_{\bm i}^{(l)}$ according to (\ref{opt:hatc})
\STATE Solve the minimization problem (\ref{opt:leastSQR}) for $\hat{\bm \alpha}_{n+1}$ and $\hat{\bm \beta}_{n+1}$
\STATE If necessary, apply the \emph{ah hoc} correction of Remark \ref{rem:AdHoc}
\STATE Let $\T_{n+1} = \T_{n} \bigcup \T^{\hat{\bm \alpha}_{n+1}, \hat{\bm \beta}_{n+1}}(L_{n+1})$ according to (\ref{sg:talphabeta})
\STATE If necessary, find $L_{n+1}$ large enough to exploit parallelism
\ENDFOR
\end{algorithmic}
\end{algorithm}

Evaluating $f_{\Lambda_n}(\bm y)$ is computationally cheap and the coefficients $\hat c_{\bm \nu}$ are computed using only $f_{\Lambda_n}(\bm y)$, however, for a high dimensional problem this process can still take $10$-$20$ minutes on $6$-core CPU. A potentially cheaper alternative is to represent the interpolant using Newton hierarchical polynomials. We define the hierarchical basis
\begin{equation}
\label{eq:1Dhbasis}
	h_1(y) = 1, \quad \mbox{and for $j>1$,} \quad h_j(y) = \prod_{i=1}^{j-1} \frac{y - y_i}{y_j - y_i}, \qquad
	H_{\bm j}(\bm y) = \prod_{k=1}^d h_{j_k}
\end{equation}
Then, for any lower tensor set $\T(L)$ and $\T_m(L)$ as defined in (\ref{sg:nodes}), there are surplus coefficients $\{ s_{\bm j} \}_{\bm j \in \T_m(L) }$ satisfying the system of equations $\sum_{ \bm 1 \leq \bm j \leq \bm i } s_{\bm j}  H_{\bm j}(\bm y_{\bm i}) = f(\bm y_{\bm i})$ for every $\bm i \in \T_m(L)$ and the interpolant can be written as
\begin{equation}
\label{sg:surpInterp}
	I_{\T(L)}^m[f](\bm y) = \sum_{ \bm j \in \T_m(L) } s_{\bm j} H_{\bm j}(\bm y)
\end{equation}
The surplus coefficients are much cheaper to compute than the Legendre coefficients, and $s_{\bm j}$ can be used as an alternative to $\hat c_{\bm \nu}$. When the one dimensional growth functions is $m(l) = l + 1$, then we can infer $\hat{\bm \alpha}$ and $\hat{\bm \beta}$ from \begin{equation}
\label{sg:lsqSurp}
	\min_{ \bm \alpha, \bm \beta, \hat C } \frac{1}{2} \sum_{\bm j \in \T_m(L)} \Big( \hat C + \bm \alpha \cdot \bm j + \bm \beta \cdot \log( \bm j + \bm 1 ) + \log( | s_{\bm j} | ) \Big)^2
\end{equation}
i.e., using $s_{\bm j}$ in place of $\hat c_{\bm \nu}$ in (\ref{opt:leastSQR}). However, this approach is not suited for the case when $m(l) > l + 1$. Any interpolant  can be written in surplus form, and the sub-ordering of $\{ y_j \}_{j = m(l-1)+1}^{m(l)}$ does not affect the nodes, polynomial space, or Lebesgue constant. However, the surpluses do depend on the sub-ordering, and therefore, $s_{\bm j}$ are influenced by $\lambda_l$ associated with growth $m(l) = l + 1$ rather than $m(l) > l+1$.

\begin{rem}\label{rem:whyLower}
When $\T(L)$ is a lower index set the three constructions (\ref{sg:mDGenInterp}), (\ref{sg:explicitI}) and (\ref{sg:surpInterp}) result in identical interpolants. In addition, the associated polynomial space is determined solely by $ \T(L)$ and the growth sequence $m(l)$, i.e., the choice of nodes does not affects the range of the operator (only the norm). The three formulas can be generalized to not lower tensor sets $ \T(L)$, however, if $ \T(L)$ is not a lower set, then the three approximations are not equal and only (\ref{sg:surpInterp}) gives an interpolant. Furthermore, when $ \T(L)$ is not lower, the associated polynomial spaces depend on the specific choice of $y_j$. In our experiments, interpolants constructed from not lower tensor sets are less accurate and we restrict our attention to lower sets $ \T(L)$.
\end{rem}

\section{Numerical results} \label{sec:numerical}

In this section we present several numerical examples using functions $f(\bm y)$ that depend on the solutions to discretized linear and nonlinear parametric PDEs. We compare the convergence rates of SG interpolants for different selections of the tensor sets, including total degree space and the $\T_n$ constructed using  Algorithm \ref{alg:ADAMA}. We also test the performance of the one dimensional rules from Table \ref{table:rules}.

\subsection{Parametrized elliptic equation}

The first three examples use similar setup involving the parametrized elliptic equation defined by
\begin{equation}
\label{num:PoissonPDE}
	\left\{  \begin{array}{rll}
		 - \nabla_{\bm x} \cdot \left( a( \bm x, \bm y ) \nabla_{\bm x} u( \bm x, \bm y ) \right) & = b(\bm x), & \bm x \in D \\
 		 u(\bm x, \bm y ) & = 0, &  \bm x \in \partial D,
	\end{array} \right.
\end{equation}
where $ \bm y \in \Gamma $ and $D = [0,1] \otimes [0,1]$. The parametrized coefficient $a(\bm x, \bm y)$ is such that
\begin{equation}
\nonumber
	0 < \min_{\bm x, \bm y} a(\bm x, \bm y) \leq \max_{\bm x, \bm y} a(\bm x, \bm y) < \infty,
\end{equation}
in which case for every $\bm y \in  \Gamma$ exists a unique $u(\bm x, \bm y) \in \mathbb{H}^1_0( D )$ that satisfies (\ref{num:PoissonPDE}), e.g., \cite{nobile2008sparse}. For a bounded functional $Q : \mathbb{H}^1_0( D ) \to \mathbb{R}$ we define the quantity of interest (QoI)
\begin{equation}
\nonumber
	f( \bm y ) = Q\left( u(\bm x, \bm y ) \right) :  \Gamma \to \mathbb{R}.
\end{equation}
The first three numerical examples differ only in the choice of $ a( \bm x, \bm y ) $, $b(\bm x)$ and $ Q\left( u(\bm x, \bm y ) \right)$. The above setup has been the thoroughly studies in literature, e.g., \cite{nobile2008sparse,Beck2014732,chkifa2014high,nobile2008anisotropic,cohen2010convergence,cohen2011analytic}, and makes for a good test bed of novel techniques.

%%%%%%%%%%%%%%%%%%%%%%%%%%%%%%%%%%%%%%%%%%%%%%%%%%%%%%%%%%%%%%%%%%%%%%%%%%%%%%%%%%%%%%%%%%%%%%%%%%%%%%%%%%%%%%%%%%%%%%%%%%%%
%%%%%%%%%%%%%%%%%%%%%%%%%%%%%%%%%%%%%%%%%%%%%%%%%%%%%%%%%%%%%%%%%%%%%%%%%%%%%%%%%%%%%%%%%%%%%%%%%%%%%%%%%%%%%%%%%%%%%%%%%%%%
%
%
%   SECTION: K-L problem
%
%
%%%%%%%%%%%%%%%%%%%%%%%%%%%%%%%%%%%%%%%%%%%%%%%%%%%%%%%%%%%%%%%%%%%%%%%%%%%%%%%%%%%%%%%%%%%%%%%%%%%%%%%%%%%%%%%%%%%%%%%%%%%%
%%%%%%%%%%%%%%%%%%%%%%%%%%%%%%%%%%%%%%%%%%%%%%%%%%%%%%%%%%%%%%%%%%%%%%%%%%%%%%%%%%%%%%%%%%%%%%%%%%%%%%%%%%%%%%%%%%%%%%%%%%%%

\subsubsection{Karhunen-Lo\'eve expansion}

Let $\left( \Omega, \mathcal{F}, \mathbb{P} \right)$ denote a complete probability space, with sample space $\Omega$, $\sigma$-algebra $\mathcal{F} = 2^\Omega $ and probability measure $\mathbb{P} : \mathcal{F} \to [0,1]$. For $x \in [0,1]$ and $\eta \in \Omega$. Let $a_1( x, \eta )$ define a random field with covariance function
\begin{equation}
\nonumber
	\mbox{Cov} \Big[ \log( a_1 - 0.5 ) \Big] ( x, x' ) = \exp\left( - \frac{(x - x')^2}{4} \right)
\end{equation}
Using Karhunen-Lo\'eve expansion, we parametrize the random field using the dominant seven eigenvalues and eigenfunctions
\begin{equation}
\nonumber
	a_1( x, \bm y) \approx 0.5 + \exp \left( 1 + \frac{\sqrt[4]{9\pi}}{2} y_1 + \frac{\sqrt[4]{9\pi}}{\sqrt{2}} \sum_{k=1}^3 e^{ -\frac{(k\pi)^2}{32} } \big( y_{2k} \sin( k \pi x_1 ) + y_{2k+1} \cos( k \pi x_1 ) \big)  \right),
\end{equation}
where $\bm y$ is uniformly distributed over $\Gamma$\cite{nobile2008sparse,nobile2008anisotropic}. Assuming that the diffusion coefficient in (\ref{num:PoissonPDE}) depends only on the first component of $\bm x$, we define $a_1(\bm x, \bm y) = a_1( x_1, \bm y)$. The source term is $b_1(\bm x) = \cos( x_1 ) \sin( x_2 )$ and the quantity of interest is the ${L}^2(D)$ norm of the solution $u(\bm x, \bm y)$
\begin{equation}
\label{num:QoIl2norm}
	f_1( \bm y ) = \left(  \int_{ D } u^2( \bm x, \bm y ) d\bm x \right)^{1/2}, \qquad \bm y \in \Gamma \subset \mathbb{R}^7
\end{equation}
Analytic expression for $f_1(\bm y)$ is not available, however, for a given $\bm y_{\bm j} \in \Gamma$, we can approximate (\ref{num:PoissonPDE}) with a finite element method (e.g., \cite{strang1973analysis}) and compute a sufficiently accurate approximation to $f_1( \bm y_{\bm j} )$. The error associated with this additional approximation is beyond the scope of this paper and we focus our attention to the discretized $f_1(\bm y)$.

\begin{figure}
\includegraphics[scale=0.5]{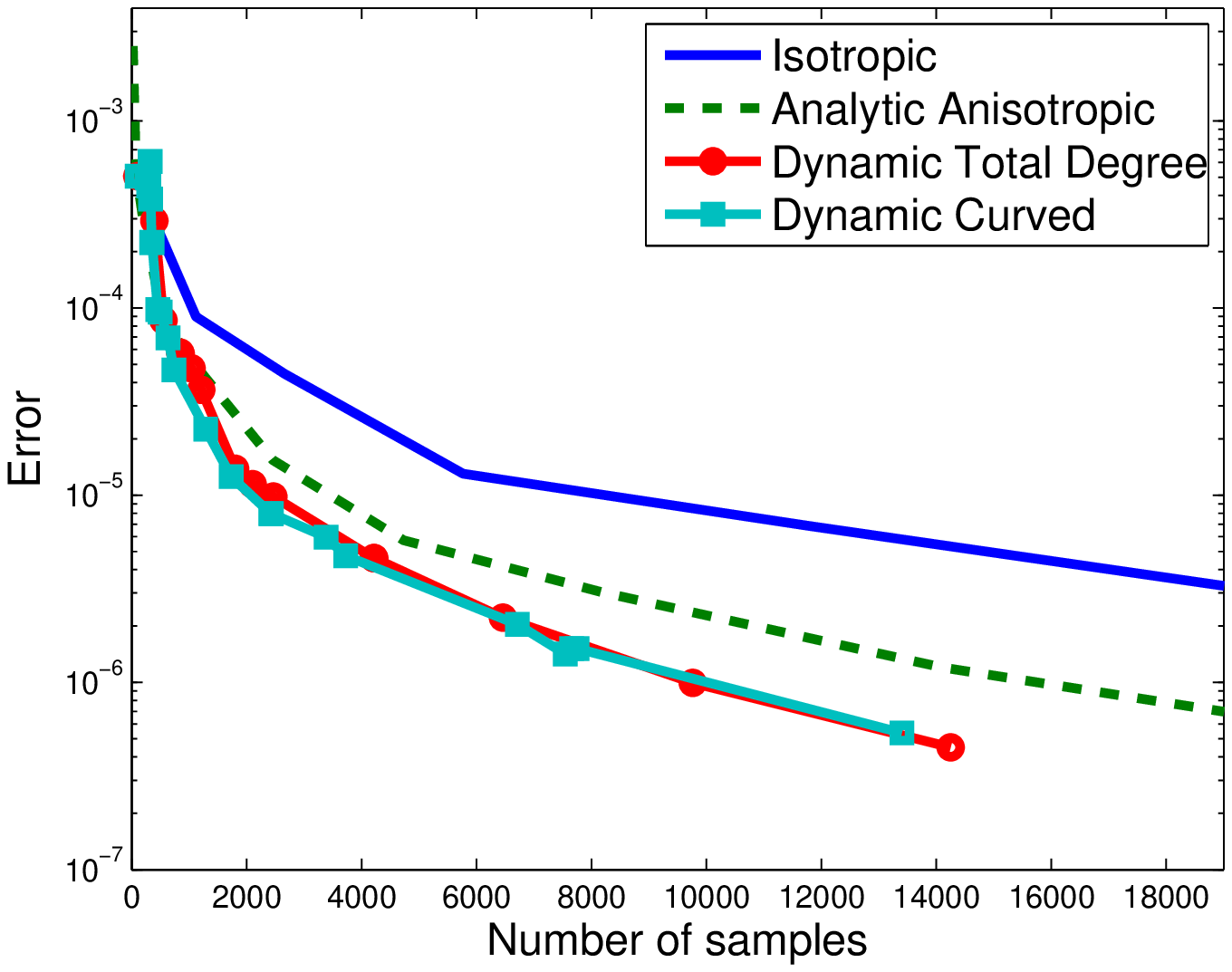}
\includegraphics[scale=0.5]{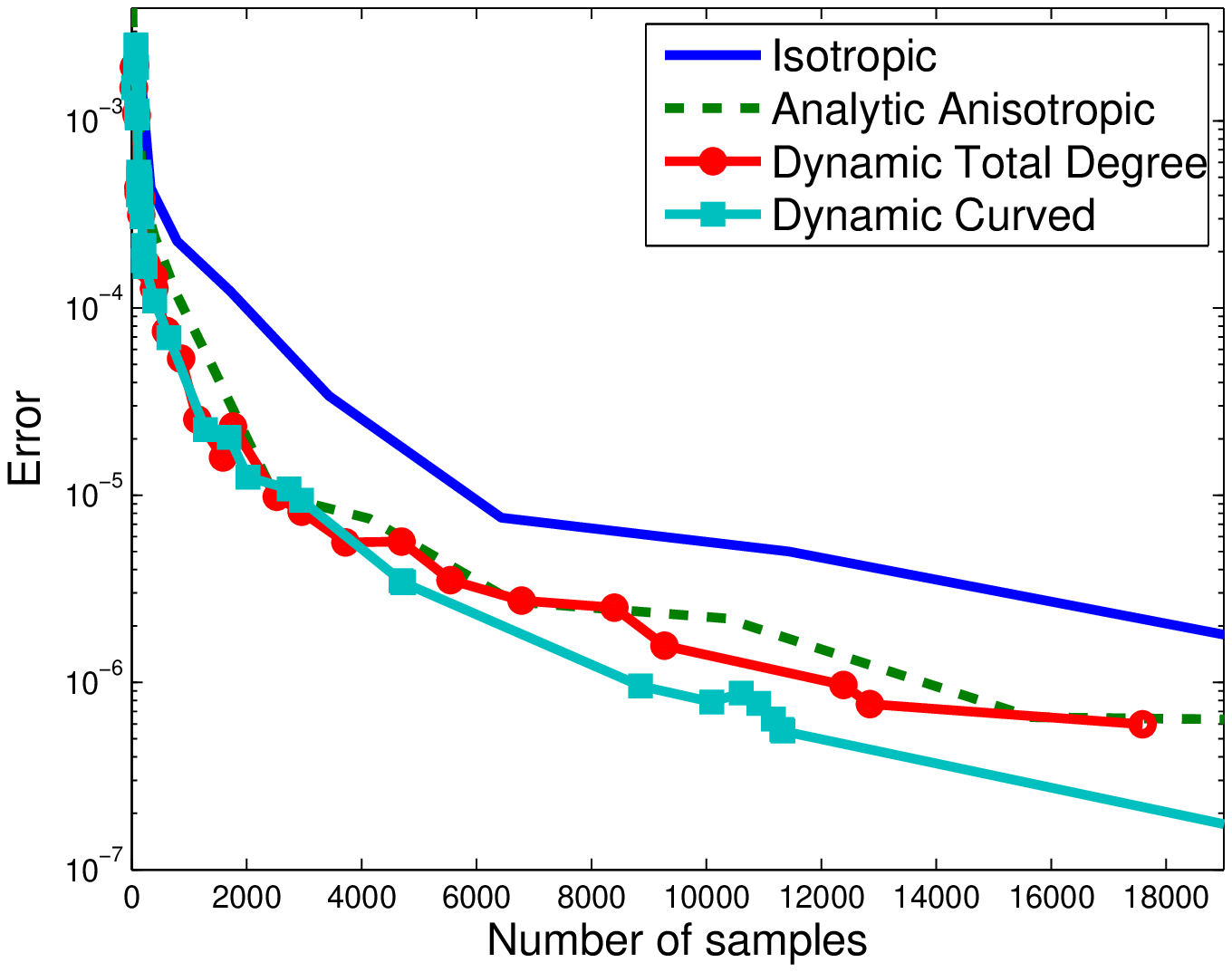}
\caption{Applying four different construction of interpolants for $f_1(\bm y)$. Left: using Clenshaw-Curtis nodes. Right: using Leja nodes.}\label{fig:KL7D_selection}
\end{figure}

Figure \ref{fig:KL7D_selection} shows a comparison of eight different types of interpolation methods applied to $f_1(\bm y)$. The left plot shows the results using the Clenshaw-Curtis rule, while the right plot uses the Leja nodes. Each plot shows the results from four different selections of the tensor set and the error is estimated from $1000$ uniformly distributed random samples $\{ \tilde{\bm y}_i \}_{i=1}^{1000}$
\begin{equation}
\label{num:mcLinf}
	error = \max_{ 1 \leq i \leq 1000 } \left| I_{ \T_n}[f_1](\tilde{\bm y}_i) - f_1(\tilde{\bm y}_i) \right|.
\end{equation}

The isotropic case corresponds to the construction of isotropic total degree polynomial space, i.e., as defined in Remark \ref{rem:growth}. However, the $e^{-k^2}$ term in (\ref{num:aKL}) decays fast with increasing $k$ and thus $y_k$ variables corresponding to larger $k$ will have lesser effect on the overall variation of $f_1(\bm y)$. The isotropic approach fails to capture the different behavior of $y_k$ and hence it is the worst performing scheme for this example. Anisotropic SG method has been proposed in \cite{nobile2008anisotropic} and using $a_1(x,\bm y)$ anisotropic weights are derived
\begin{equation}
\nonumber
	\bm \alpha = \left( 0.85, 0.80, 0.80, 1.0, 1.0, 1.6, 1.6 \right),
\end{equation}
and the analytic anisotropic polynomial space is defined by
\begin{equation}
\nonumber
	{\bm \Lambda}^{\bm \alpha}(L) = \big\{ \bm \nu \in \mathbb{N}^d : {\bm \alpha} \cdot \bm \nu \leq L \big\}.
\end{equation}
Using ${\bm \Lambda}^{\bm \alpha}(L)$, we construct the corresponding interpolant defined in Theorem \ref{theorem:opt}, the results are plotted as the dashed lines in Figure \ref{fig:KL7D_selection}.

For comparison purposes, Figure \ref{fig:KL7D_selection} also plots a \emph{Dynamic Total Degree} example that uses a modified version of Algorithm \ref{alg:ADAMA}. In the total degree case, we remove the $\bm \beta$ term from the weight function (\ref{opt:decayModel}) and the least-squares estimate (\ref{opt:leastSQR}), and we compute only $\hat{\bm \alpha}$. For both Clenshaw-Curtis and Leja nodes, the total degree approach captures the anisotropic behavior of $f_1(\bm y)$ and the method exhibits faster convergence. The estimated normalized parameters $\hat{\bm \alpha}$ are given on Table \ref{table:KLalphaTD}. Since the dominant polynomial space is independent from the scaling of $\hat{\bm \alpha}$, we divide all entries of $\hat{\bm \alpha}$ by the average $\frac{1}{7}|\hat{\bm \alpha}|_1$. Both Clenshaw-Curtis and Leja nodes result in nearly identical $\hat{\bm \alpha}$ even though Clenshaw-Curtis parameters are based on the projection of the interpolant (\ref{opt:hatc}) while the Leja parameters are computed using the hierarchical surpluses (\ref{sg:lsqSurp}). The Clenshaw-Curtis nodes have a lower Lebesgue constant, however, the Clenshaw-Curtis nodes result in grids with more nodes due to the exponential growth of $m(l)$. For this particular problem using dynamic total degree nodes, the lower operator norm of Clenshaw-Curtis nodes leads to faster convergence.

Figure \ref{fig:KL7D_selection} also shows the results from applying Algorithm \ref{alg:ADAMA} to $f_1(\bm y)$, and the curve is labeled \emph{Dynamic Curved}. When using Clenshaw-Curtis nodes, there is no significant difference between the total degree method and the quasi-optimal interpolation (\ref{opt:optInterp}). The estimated (normalized) decay parameters $\hat{\bm \alpha}$ and $\hat{\bm \beta}$ are listed in Table \ref{table:KLbeta} and we see that $\hat{\bm \beta}$ is relatively small. In contrast, when using the Leja nodes, the quasi-optimal estimate (\ref{opt:optInterp}) leads to a significant improvement in convergence and from Table \ref{table:KLbeta} we see that the $\hat{\bm \beta}$ estimated parameters are smaller compared to the Clenshaw-Curtis ones, which is due to the larger Lebesgue constant.

\begin{table}
\begin{center}
\begin{tabular}{|l|r|r|r|r|r|r|r|}
\hline
Dimension & $\hat \alpha_1$ & $\hat \alpha_2$ & $\hat \alpha_3$ & $\hat \alpha_4$ & $\hat \alpha_5$ & $\hat \alpha_6$ & $\hat \alpha_7$ \\
\hline
Clenshaw-Curtis & $0.58$ & $0.70$ & $0.57$ & $0.99$ & $0.93$ & $1.57$ & $1.66$ \\
Leja            & $0.61$ & $0.70$ & $0.58$ & $0.96$ & $0.91$ & $1.55$ & $1.68$ \\
Analytic        & $0.78$ & $0.73$ & $0.73$ & $0.92$ & $0.92$ & $1.46$ & $1.46$ \\
%Analytic        & ${\color{red} 0.78}$ & $0.73$ & $0.73$ & $0.92$ & $0.92$ & $1.46$ & $1.46$ \\
\hline
\end{tabular}
\caption{Estimated $\hat{\bm \alpha}$ parameters for $f_1(\bm y)$ using total degree estimate.}\label{table:KLalphaTD}
\end{center}
\end{table}

\begin{table}
\begin{center}
\begin{tabular}{|l|r|r|r|r|r|r|r|}
\hline
Dimension & $\hat \alpha_1$ & $\hat \alpha_2$ & $\hat \alpha_3$ & $\hat \alpha_4$ & $\hat \alpha_5$ & $\hat \alpha_6$ & $\hat \alpha_7$ \\
\hline
Clenshaw-Curtis & $0.64$ & $0.77$ & $0.80$ & $0.97$ & $0.99$ & $1.48$ & $1.34$ \\
Leja            & $0.67$ & $0.80$ & $0.85$ & $1.08$ & $1.05$ & $1.43$ & $1.12$ \\
\hline
Dimension & $\hat \beta_1$ & $\hat \beta_2$ & $\hat \beta_3$ & $\hat \beta_4$ & $\hat \beta_5$ & $\hat \beta_6$ & $\hat \beta_7$ \\
\hline
Clenshaw-Curtis & $-0.46$ & $-0.51$ & $-0.94$ & $-0.31$ & $-0.47$ & $-0.29$ & $0.09$ \\
Leja            & $-1.00$ & $-1.14$ & $-1.68$ & $-1.24$ & $-1.27$ & $-0.96$ & $-0.21$ \\
\hline
\end{tabular}
\caption{Estimated $\hat{\bm \alpha}$ and $\hat{\bm \beta}$ parameters for $f_1(\bm y)$.}\label{table:KLbeta}
\end{center}
\end{table}

\begin{figure}
\begin{center}
\includegraphics[scale=0.75]{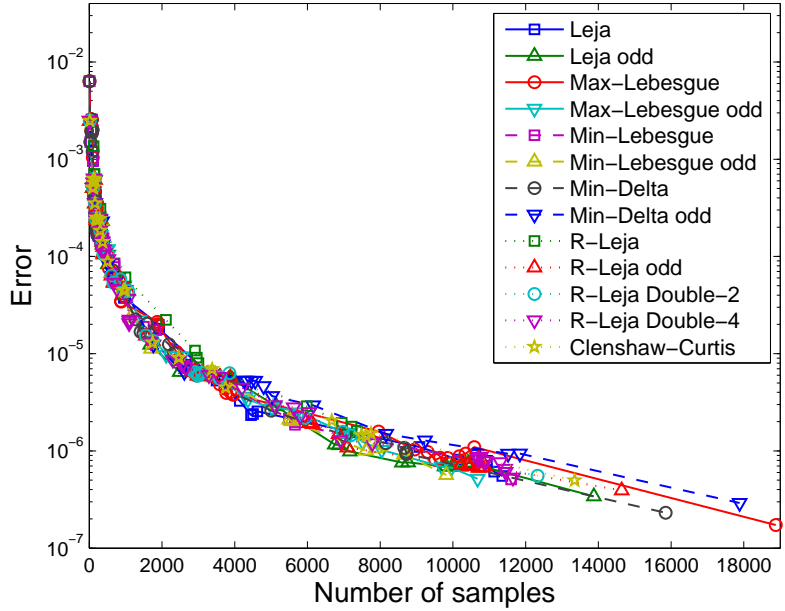}
\end{center}
\caption{Algorithm \ref{alg:ADAMA} applied to $f_1(\bm y)$ using all rules from Table \ref{table:rules}.}\label{fig:KL7D_allRules}
\end{figure}

%\begin{figure}
%\includegraphics[scale=0.5]{figs/allAlpha.eps}
%\includegraphics[scale=0.5]{figs/allBeta.eps}
%\caption{The variation in the values of the estimated $\hat{\bm \alpha}$ and $\hat{\bm \beta}$ using interpolation rules with different $\lambda_l$.}\label{fig:KL7D_allAB}
%\end{figure}

\begin{table}
\begin{center}
\begin{tabular}{|l|r|r|r|r|r|r|r|r|r|r|}
\hline
Dimension & $\delta \hat \alpha_1$ & $\delta \hat \alpha_2$ & $\delta \hat \alpha_3$ & $ \delta \hat \alpha_4$ & $ \delta \hat \alpha_5$ & $\delta \hat \beta_1$ & $\delta \hat \beta_2$ & $\delta \hat \beta_3$ & $\delta \hat \beta_4$ & $\delta \hat \beta_5$  \\
\hline
Clenshaw-Curtis & $0.80$ & $0.45$ & $0.82$ & $0.50$ & $0.67$ & $1.27$ & $0.88$ & $0.93$ & $1.70$ & $0.89$ \\
%Leja            & $0.67$ & $0.80$ & $0.85$ & $1.08$ & $1.05$ & $1.43$ & $1.12$ \\
%\hline
%Dimension &  & $\hat \beta_6$ & $\hat \beta_7$ \\
%\hline
%Clenshaw-Curtis & $-0.46$ & $-0.51$ & $-0.94$ & $-0.31$ & $-0.47$ & $-0.29$ & $0.09$ \\
%Leja            & $-1.00$ & $-1.14$ & $-1.68$ & $-1.24$ & $-1.27$ & $-0.96$ & $-0.21$ \\
\hline
\end{tabular}
\caption{Relative variation of $\hat{\bm \alpha}$ and $\hat{\bm \beta}$ parameters for $f_1(\bm y)$ over all rules from  Table \ref{table:rules} .}\label{table:AllAB}
\end{center}
\end{table}

Figure \ref{fig:KL7D_allRules} shows the results of applying Algorithm \ref{alg:ADAMA} to $f_1(\bm y)$ and using different rules from Table \ref{table:rules}. The curves are not smooth since the results are affected by error in the least-squares approximation, fluctuations of the operator norm (see Figure \ref{fig:LebPenalty}), the random samples used to compute the error (\ref{num:mcLinf}), and, when $m(l) > l+1$ the range of the resulting interpolant is a super-set of the estimated optimal polynomial space. However, all interpolants have similar convergence rate within the ``noise'' in the method. Table \ref{table:AllAB} gives the relative variation of $\hat{\bm \alpha}$ and $\hat{\bm \beta}$ components for all rules in Table \ref{table:rules}, where $\delta \alpha_k$ is computed as the difference between the largest and smallest estimate for $\alpha_k$ divided by the average (same for $\beta_k$). From the first $5$ components, we see that the components of $\hat{\bm \beta}$ vary over a much wider range than those of $\hat{\bm \alpha}$, since the $\hat{\bm \alpha}$ is affected by only by the ``noise'' while $\hat{\bm \beta}$ is also affected by the different Lebesgue constants for the different rules. Note that the table lists only the first $5$ components, since $y_6$ and $y_7$ are associated with the smallest eigenvalue of the Karhunen-Lo\'eve expansion, $f_1(\bm y)$ is least sensitive to $y_6$ and $y_7$, and the constructed interpolants have few nodes in those two directions, which leads to unreliable estimates.

%{\bf Fix this} Figure \ref{fig:KL7D_allAB} plots the fluctuation of the estimates $\hat{\bm \alpha}$ and $\hat{\bm \beta}$ for the different interpolation rules, while the $\hat{\bm \alpha}$ parameters are fairly consistent, the $\hat{\bm \beta}$ parameters exhibit larger fluctuation corresponding to the difference in $\lambda_l$. Note, the largest fluctuation corresponds to inputs $y_6$ and $y_7$, which are associated with the smallest eigenvalue of the Karhunen-Lo\'eve expansion. The $f^{(7)}(\bm y)$ is least sensitive to $y_6$ and $y_7$, and the constructed interpolants have few nodes in those two directions, hence, the least squares error in estimating the corresponding parameters is the largest.

%%%%%%%%%%%%%%%%%%%%%%%%%%%%%%%%%%%%%%%%%%%%%%%%%%%%%%%%%%%%%%%%%%%%%%%%%%%%%%%%%%%%%%%%%%%%%%%%%%%%%%%%%%%%%%%%%%%%%%%%%%%%
%%%%%%%%%%%%%%%%%%%%%%%%%%%%%%%%%%%%%%%%%%%%%%%%%%%%%%%%%%%%%%%%%%%%%%%%%%%%%%%%%%%%%%%%%%%%%%%%%%%%%%%%%%%%%%%%%%%%%%%%%%%%
%
%
%   SECTION: Inclusion problem
%
%
%%%%%%%%%%%%%%%%%%%%%%%%%%%%%%%%%%%%%%%%%%%%%%%%%%%%%%%%%%%%%%%%%%%%%%%%%%%%%%%%%%%%%%%%%%%%%%%%%%%%%%%%%%%%%%%%%%%%%%%%%%%%
%%%%%%%%%%%%%%%%%%%%%%%%%%%%%%%%%%%%%%%%%%%%%%%%%%%%%%%%%%%%%%%%%%%%%%%%%%%%%%%%%%%%%%%%%%%%%%%%%%%%%%%%%%%%%%%%%%%%%%%%%%%%

\subsubsection{Inclusion problem}

\begin{figure}
\begin{center}
\includegraphics[scale=0.5]{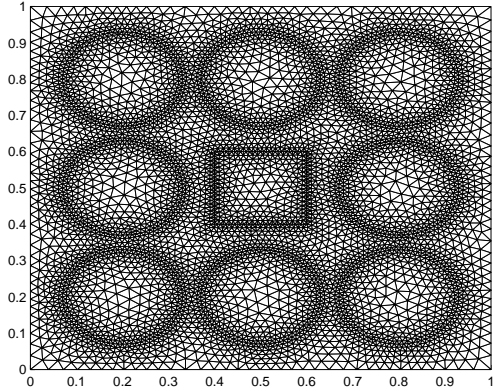}
\end{center}
\caption{Mesh for the inclusion problem. Each disk has radius $\frac{2}{15}$ and the centers are at $0.2$, $0.5$ and $0.8$ along the horizontal and vertical axis. The square box is located at the center and has side of $0.2$.}\label{num:diskmesh}
\end{figure}

In this section we apply our approach to two inclusion problems, presented in \cite{Beck2014732}, where we consider the domain given on Figure \ref{num:diskmesh}.

\subsection*{Isotropic case}
\begin{figure}
\includegraphics[scale=0.5]{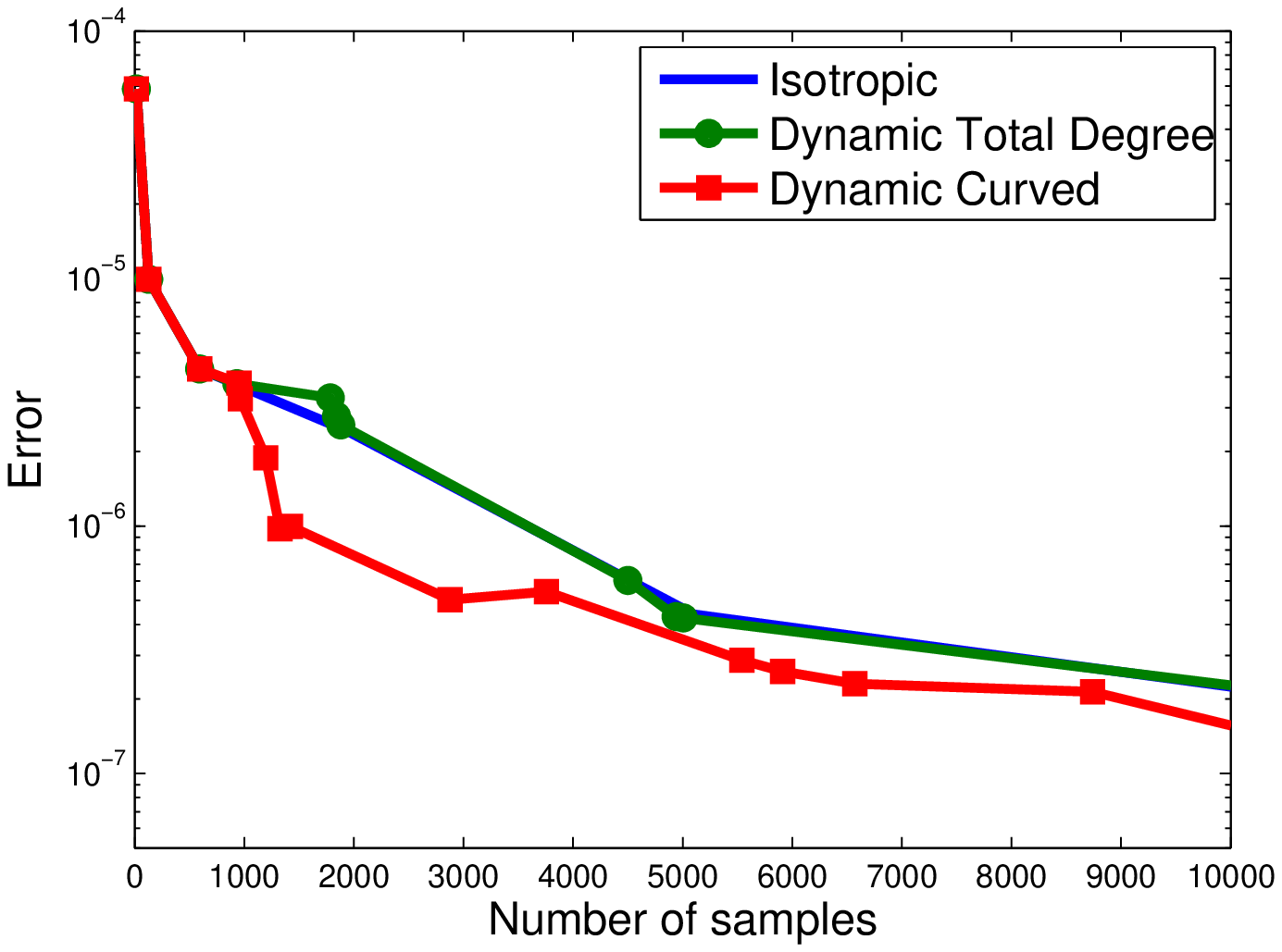}
\includegraphics[scale=0.5]{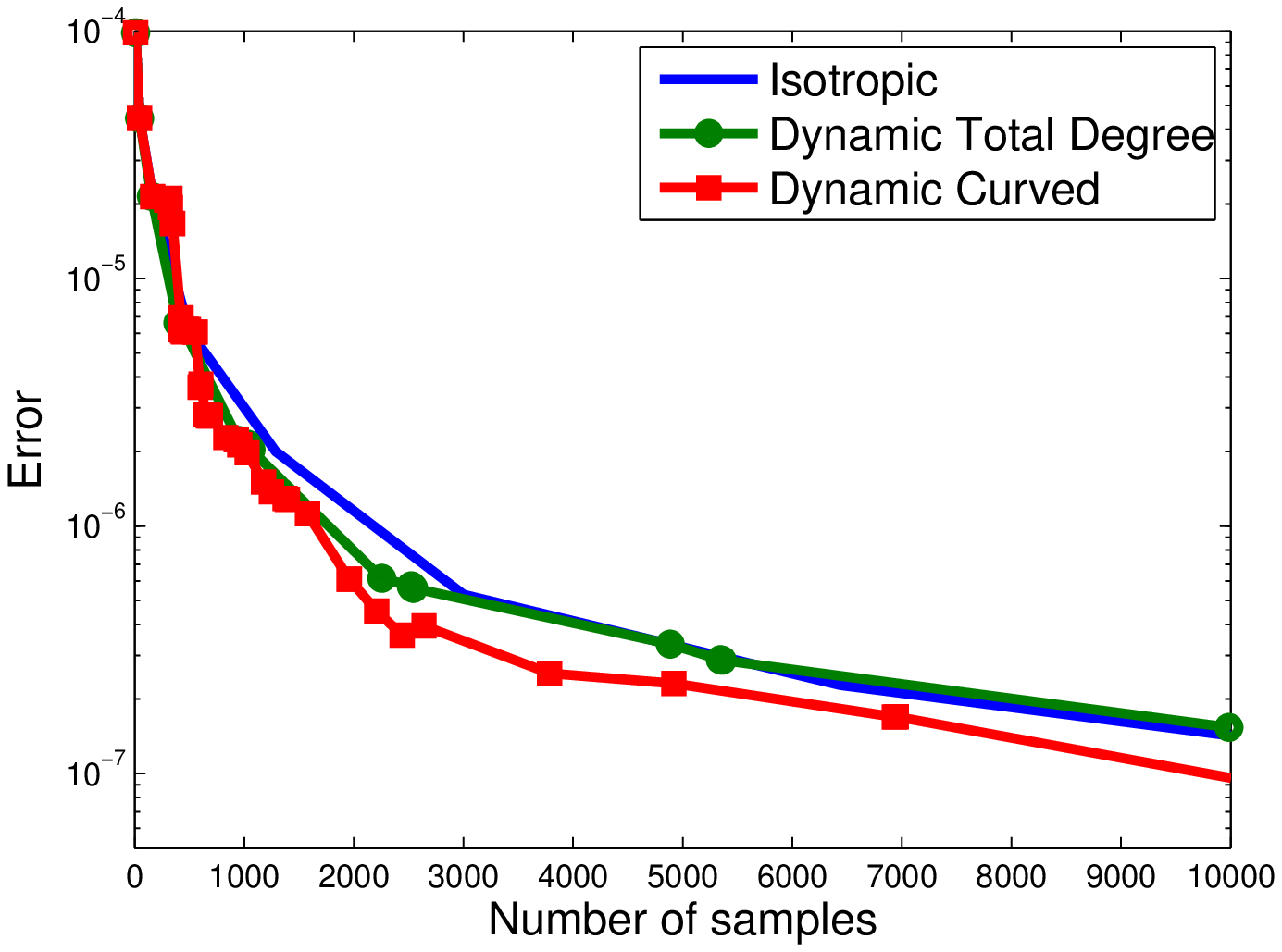}
\caption{Applying three different construction of interpolants for $f_2(\bm y)$. Left: using Clenshaw-Curtis nodes. Right: using Leja nodes.}\label{fig:IN8D_selection}
\end{figure}

The (almost) isotropic inclusion problem associates each of the eight disks with a component of $\bm y \in \bigotimes_{k=1}^8 [0.01,1]$ and we define
\begin{equation}
\nonumber
	a_2( \bm x, \bm y ) = \left\{ \begin{array}{ll}
		y_k, & \bm x \in \mbox{Disk} \, k \\
		1, & \mbox{otherwise}.
	\end{array} \right.
	\quad
	b_2(\bm x) = \left\{ \begin{array}{ll}
		100, & \bm x \in \mbox{Box} \\
		0, & \mbox{otherwise}.
	\end{array} \right.
	\quad
	f_2(\bm y) = \int_{\mbox{Box}} u(\bm x, \bm y) d\bm x
\end{equation}
Figure \ref{fig:IN8D_selection} shows a plot of the convergence rate of interpolants based on Clenshaw-Curtis and Leja nodes. Figure \ref{fig:IN8D_selection} compares three possible selections of the polynomial space, isotropic total degree (see Remark \ref{rem:growth}), dynamic curved (i.e., Algorithm \ref{alg:ADAMA}), and dynamic total degree (i.e., Algorithm \ref{alg:ADAMA} ignoring $\bm \beta$). The disks are not equidistant from the box and thus the components of $\bm y$ do not have equal influence on $f_2(\bm y)$, however, the anisotropy is weak and the performance of the the anisotropic total selection is indistinguishable from the isotropic interpolant. On the other hand, our quasi-optimal approach has a significantly better convergence. Table \ref{table:QoI8disk} gives the final values of the estimated $\hat{\bm \alpha}$ and $\hat{\bm \beta}$ parameters, similar to the Karhunen-Lo\'eve example, the $\hat{\bm \alpha}$ parameters are almost identical for both Clenshaw-Curtis and Leja nodes, however, the $\hat{\bm \beta}$ parameters for the Leja nodes are smaller due to the larger operator norm. Despite the seemingly small values of $\hat{\bm \beta}$, the added flexibility of the $\log$ term in (\ref{opt:optInterp}) results in accurate interpolants with thousands of fewer nodes as compared to the total degree construction.

It is interesting to note that several of the $\hat \beta_k$ estimates in Table \ref{table:QoI8disk} are in fact positive, as opposed to the negative values suggested in the theoretical derivation of (\ref{opt:optInterp}) and Assumption \ref{ass:Lebesgue}. This indicates that the coefficients $\hat{c}_{\bm i}$ exhibit a combination of exponential and algebraic decay and the added flexibility of $\bm \beta$ (i.e., imposing no restictions on the individual $\hat \beta_k$), results in a more accurate estimate of the optimal polynomial space. Thus, our approach is applicable to a wider class of functions, namely those obeying estimate (\ref{opt:optInterp}) for any $\bm \beta$, albeit rigorous characterization of this class of functions is not currently available.

\begin{table}
\begin{center}
\begin{tabular}{|l|r|r|r|r|r|r|r|r|}
\hline
Dimension & $\hat \alpha_1$ & $\hat \alpha_2$ & $\hat \alpha_3$ & $\hat \alpha_4$ & $\hat \alpha_5$ & $\hat \alpha_6$ & $\hat \alpha_7$ & $\hat \alpha_8$ \\
\hline
Clenshaw-Curtis & $0.86$ & $1.32$ & $0.67$ & $1.32$ & $1.25$ & $0.67$ & $1.26$ & $0.66$ \\
Leja            & $0.94$ & $1.32$ & $0.66$ & $1.32$ & $1.22$ & $0.66$ & $1.22$ & $0.66$ \\
\hline
Dimension & $\hat \beta_1$ & $\hat \beta_2$ & $\hat \beta_3$ & $\hat \beta_4$ & $\hat \beta_5$ & $\hat \beta_6$ & $\hat \beta_7$ & $\hat \beta_8$ \\
\hline
Clenshaw-Curtis & $ 0.66$ & $-0.49$ & $ 1.41$ & $-0.50$ & $-0.20$ & $ 1.40$ & $-0.20$ & $1.47$ \\
Leja            & $ 0.17$ & $-1.24$ & $ 0.75$ & $-1.24$ & $-0.86$ & $ 0.74$ & $-0.86$ & $0.79$ \\
\hline
\end{tabular}
\caption{Estimated $\hat{\bm \alpha}$ and $\hat{\bm \beta}$ parameters for $f_2(\bm y)$.}\label{table:QoI8disk}
\end{center}
\end{table}

\subsection*{Anisotropic case}
\begin{figure}
\includegraphics[scale=0.5]{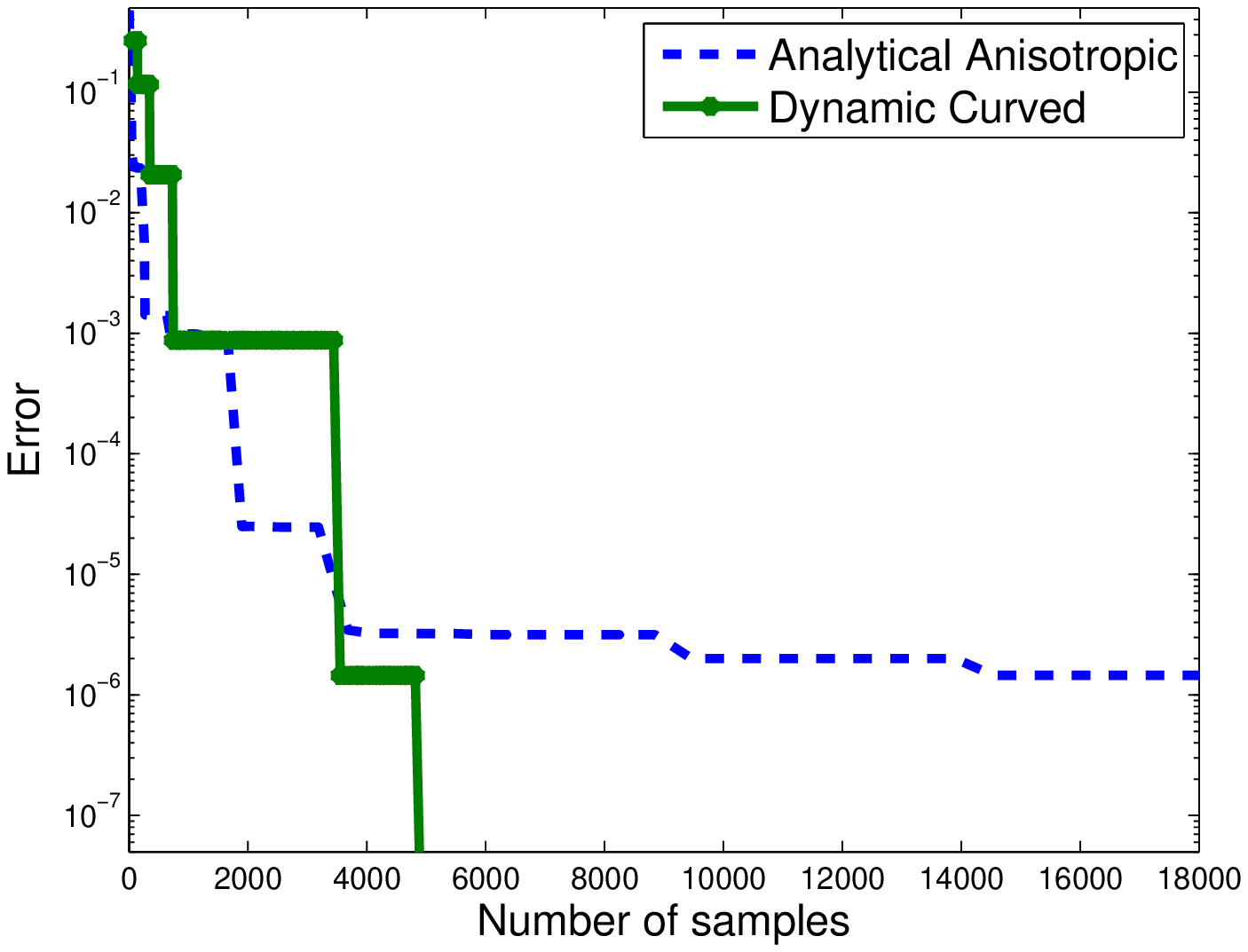}
\includegraphics[scale=0.5]{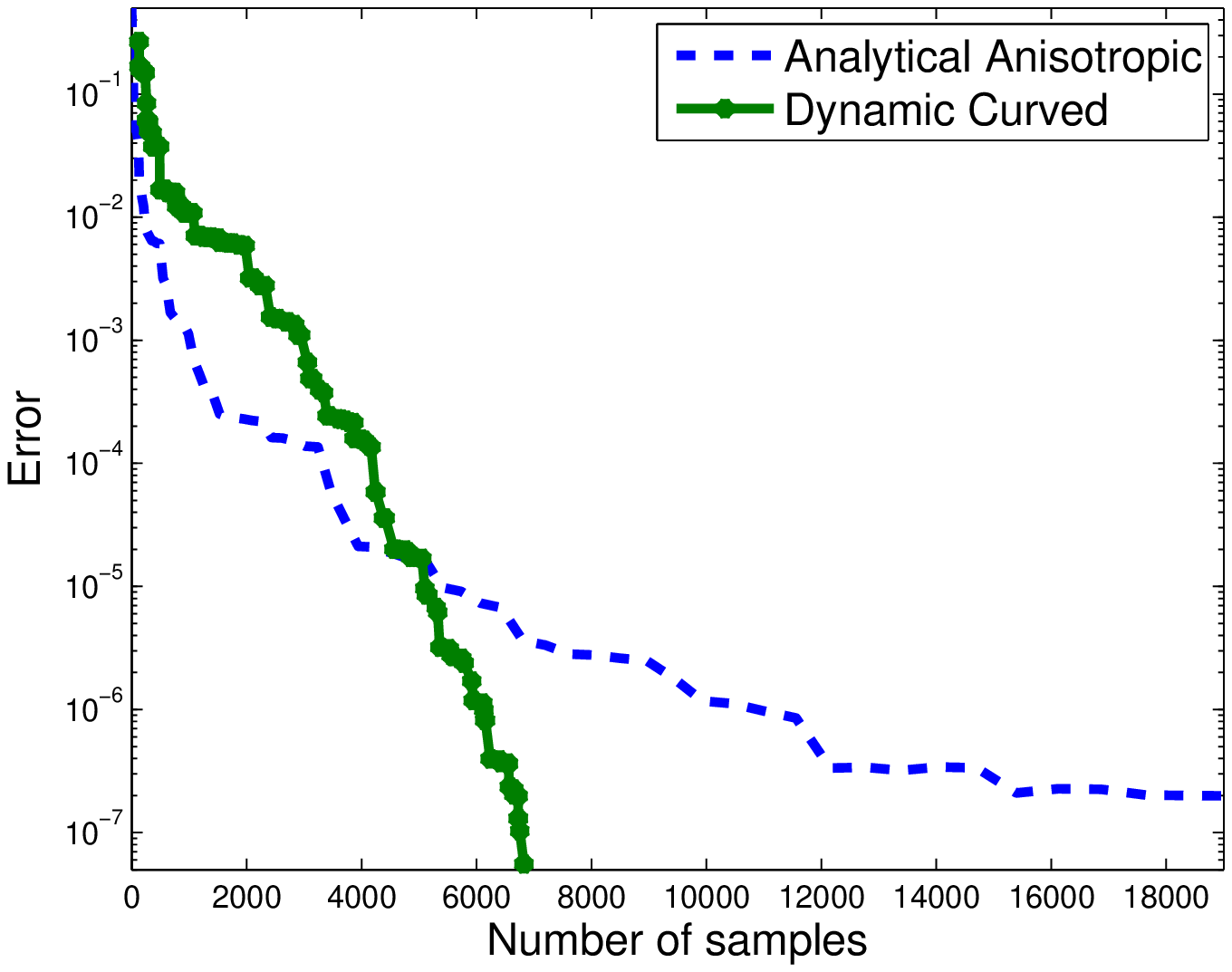}
\caption{Applying three different construction of interpolants $f_3(\bm y)$. Left: using Clenshaw-Curtis nodes. Right: using Leja nodes.}\label{fig:IN4D_selection}
\end{figure}

The second variation of the inclusion problem uses only the four corner disks
\begin{align*}
	y_{\mbox{top-left}}  \in [0.109,1], && y_{\mbox{top-right}} \in [0.2575,1], \\
	y_{\mbox{bottom-left}}  \in [0.010,1], && y_{\mbox{bottom-right}} \in [0.4060,1].
\end{align*}
The diffusion coefficient is constant $1$ outside the four boxes. The four parameter forcing term and quantity of interest are now defined over the entire domain, i.e., $b_3(\bm x) = 100$ and
\begin{equation}
\label{num:QoI4Disk}
	f_3(\bm y) = \int_{D} u(\bm x, \bm y) d\bm x.
\end{equation}
Despite the smaller number of dimensions, the four parameter problem is more difficult, due to the smaller $a_{min}$ and smaller region of analyticity. For this problem, anisotropic total degree weights have been derived
\begin{equation}
\label{num:4Dalpha}
	{\bm \alpha} = \left( 40, 317, 137, 227 \right),
\end{equation}
and numerical results show that (\ref{num:4Dalpha}) are reasonably accurate when used in the context of Galerkin projection\cite{Beck2014732}. Figure \ref{fig:IN4D_selection} shows a comparison in the convergence rate between interpolants constructed via the total degree weights given in (\ref{num:4Dalpha}) and the application of Algorithm \ref{alg:ADAMA}, Table \ref{table:QoI4isk} lists the estimated $\hat{\bm \alpha}$ and $\hat{\bm \beta}$ parameters. The convergence rate of the interpolants is closer to algebraic, leading to a large  $\bm \beta$, which dominates for indexes with small entries. Thus, in the context of interpolation, the total degree space estimated by the weights (\ref{num:4Dalpha}) is far from optimal.  Initially, the least squares method struggles to estimate the quasi-optimal parameters, however, once a sufficient number of samples are computed, the convergence rate of Algorithm \ref{alg:ADAMA} dramatically increases leading to enormous savings compared to the total degree selection. Furthermore, due to the small region of analyticity of $f_3(\bm y)$, the difference in Lebesgue constant between Clenshaw-Curtis and Leja has a significant effect with the former method outperforming the latter.

\begin{table}
\begin{center}
\begin{tabular}{|l|r|r|r|r|r|r|r|r|}
\hline
Dimension & $\hat \alpha_1$ & $\hat \alpha_2$ & $\hat \alpha_3$ & $\hat \alpha_4$ & $\hat \beta_1$ & $\hat \beta_2$ & $\hat \beta_3$ & $\hat \beta_4$ \\
\hline
Clenshaw-Curtis & $0.18$ & $2.38$ & $0.18$ & $1.27$ & $12.70$   & $10.64$   & $13.62$   & $12.00$ \\
Leja            & $0.24$ & $2.22$ & $0.24$ & $1.30$ & $\, 6.95$ & $\, 5.69$ & $\, 8.75$ & $\, 6.94$ \\
\hline
\end{tabular}
\caption{Estimated $\hat{\bm \alpha}$ and $\hat{\bm \beta}$ parameters for $f_3(\bm y)$.}\label{table:QoI4isk}
\end{center}
\end{table}

%%%%%%%%%%%%%%%%%%%%%%%%%%%%%%%%%%%%%%%%%%%%%%%%%%%%%%%%%%%%%%%%%%%%%%%%%%%%%%%%%%%%%%%%%%%%%%%%%%%%%%%%%%%%%%%%%%%%%%%%%%%%
%%%%%%%%%%%%%%%%%%%%%%%%%%%%%%%%%%%%%%%%%%%%%%%%%%%%%%%%%%%%%%%%%%%%%%%%%%%%%%%%%%%%%%%%%%%%%%%%%%%%%%%%%%%%%%%%%%%%%%%%%%%%
%
%
%   SECTION: Burgers problem
%
%
%%%%%%%%%%%%%%%%%%%%%%%%%%%%%%%%%%%%%%%%%%%%%%%%%%%%%%%%%%%%%%%%%%%%%%%%%%%%%%%%%%%%%%%%%%%%%%%%%%%%%%%%%%%%%%%%%%%%%%%%%%%%
%%%%%%%%%%%%%%%%%%%%%%%%%%%%%%%%%%%%%%%%%%%%%%%%%%%%%%%%%%%%%%%%%%%%%%%%%%%%%%%%%%%%%%%%%%%%%%%%%%%%%%%%%%%%%%%%%%%%%%%%%%%%

\subsection{Parametrized steady state Burgers equation}
\begin{figure}
\begin{center}
\includegraphics[scale=0.5]{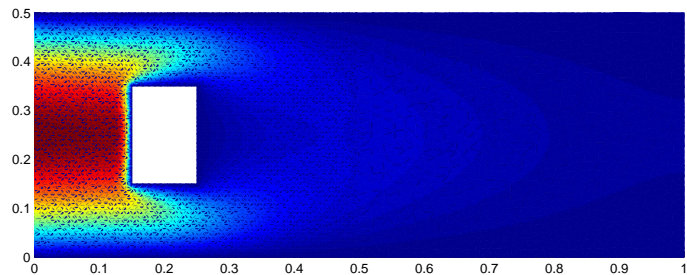}
\end{center}
\caption{The domain associated with the Burgers equation (\ref{num:burg}), the color map corresponds to the nominal solution $\bm y = \bm 0$.}\label{fig:burgmesh}
\end{figure}

Consider the steady state Burgers equation defined over the domain $D = [0,1] \otimes [0,0.5] \setminus [0.15,0.25] \otimes [0.15,0.35]$ (see Figure \ref{fig:burgmesh}). The right most wall of the domain is associated with homogenous Neumann boundary conditions, i.e.,  $\partial D_{n} = \{ 1 \} \otimes [0,0.5]$, the remainder of the boundary, i.e., $\partial D_{d}$, is associated with the Dirichlet conditions defined by
\vspace{-0.2cm}
\begin{equation}
\label{num:burgDirichlet}
	u_b(\bm x) = \left\{ \begin{array}{ll}
		16 x_2 ( 0.5 - x_2 ), & x_1 = 0 \\
		0, & \mbox{otherwise.}
	\end{array} \right.
\vspace{-0.3cm}
\end{equation}
A time dependent control problem using the above domain is described in \cite{djouadi2007optimal}. In our context, we are interested in parametrizing the steady state Burgers equation given by
\vspace{-0.2cm}
\begin{equation}
\label{num:burg}
\left\{ \begin{array}{rll}
	- \nabla_{\bm x} \cdot \left( a(\bm y) \nabla_{\bm x} u(\bm x, \bm y ) \right) + \big( \bm v  (\bm y ) \cdot \nabla_{\bm x} u(\bm x, \bm y) \big) u(\bm x, \bm y) & = 0, &  \bm x \in D, \\
	u(\bm x, \bm y) & = u_b(\bm x), & \bm x \in \partial D_d, \\
	\frac{\partial}{\partial x_1} u(\bm x, \bm y) & = 0, & \bm x \in \partial D_n,
\end{array}
\right.
\vspace{-0.3cm}
\end{equation}
where $\bm y \in \bigotimes_{k=1}^3 [-1,1]$ and
\vspace{-0.3cm}
\begin{equation}
\nonumber
	a(\bm y) = \frac{1}{200 + 100 y_3}, \qquad \bm v(\bm y) = \left( \begin{array}{r} 1 + 0.2 y_1 \\  0.1 y_2 \end{array} \right),
	\qquad
	f_4(\bm y) = \int_{0.05}^{0.45} \int_{0.6}^{0.8} u(\bm x, \bm y) dx_1 dx_2
\end{equation}
Figure \ref{fig:burgmesh} plots the nominal solution corresponding to $\bm y = \bm 0$. There is no available result regarding the analyticity of $f_4(\bm y)$ or the corresponding optimal polynomial space. However, the components of the convection term $\bm v(\bm y)$ affect different derivatives of the non-uniform solution, $y_1$ and $y_2$ are multiplied by different coefficients, and $y_3$ affects the diffusion term which has very different effect on $u(\bm x, \bm y)$. Thus, we expect anisotropic decay behavior of $f_4(\bm y)$.

\begin{figure}
\includegraphics[scale=0.5]{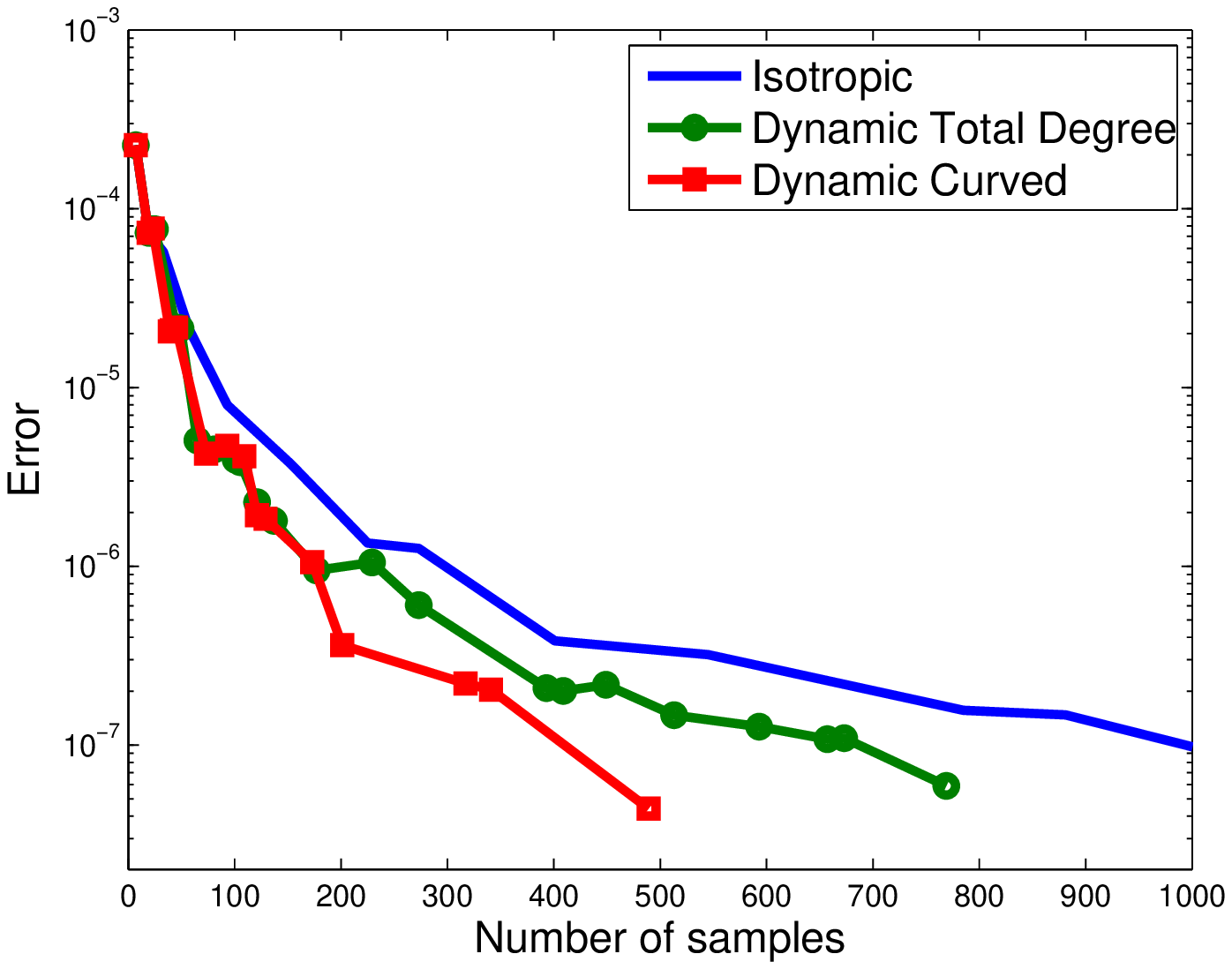}
\includegraphics[scale=0.5]{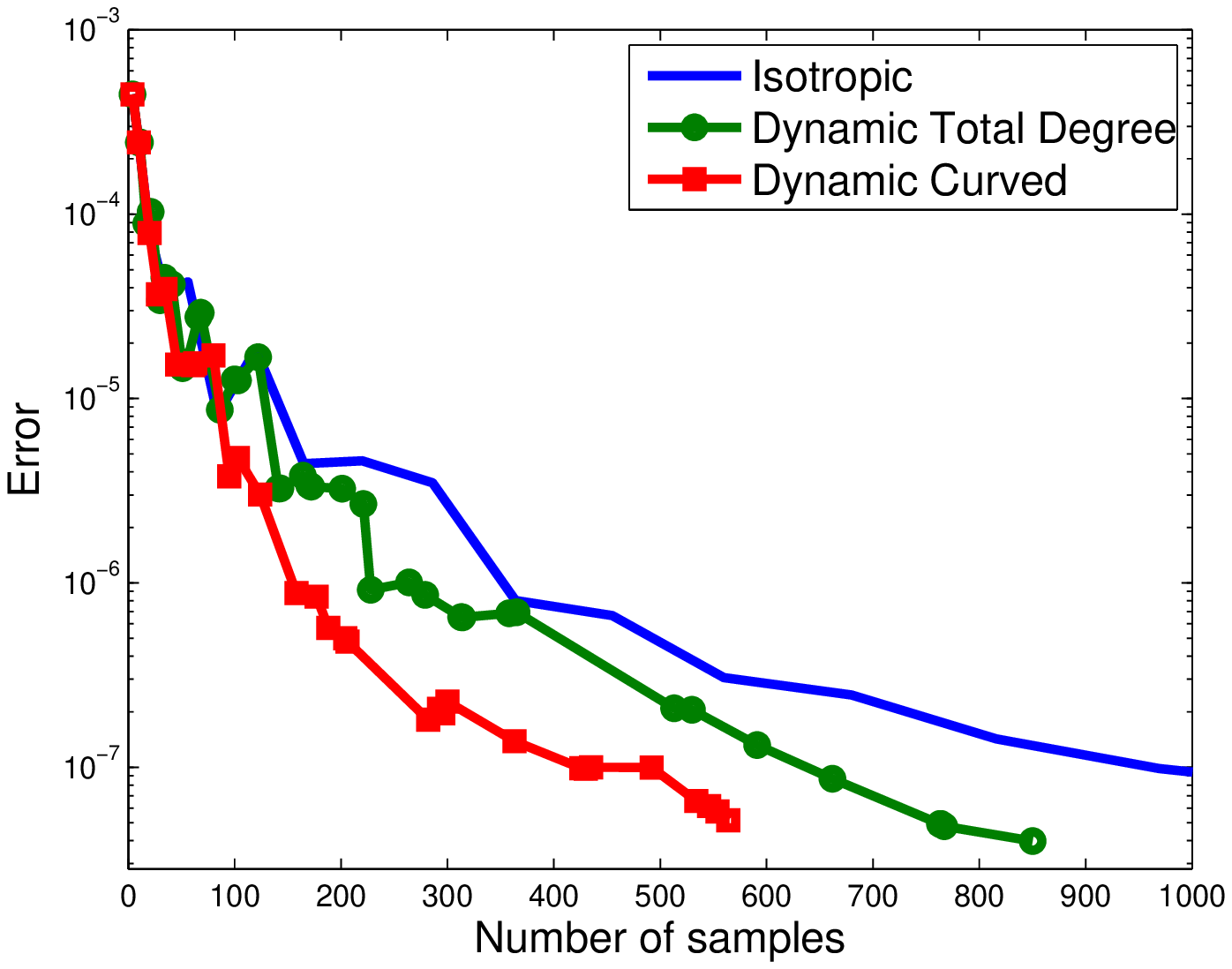}
\caption{Applying three different construction of interpolants for $f_4(\bm y)$. Left: using Clenshaw-Curtis nodes. Right: using Leja nodes.}\label{fig:burgselect}
\end{figure}

Figure \ref{fig:burgselect} plots the result of interpolating $f_4(\bm y)$ with three different polynomial selection schemes, the isotropic total degree (see Remark \ref{rem:growth}), and Algorithm \ref{alg:ADAMA} ignoring and including $\bm \beta$. We use Leja and Clenshaw-Curtis nodes and the estimated parameters are given in Table \ref{table:burg}. We observe results very similar to the previous examples, which shows that our approach extends to problems associated with nonlinear PDEs.

\begin{table}
\begin{center}
\begin{tabular}{|l|r|r|r|r|r|r|}
\hline
Dimension & $\hat \alpha_1$ & $\hat \alpha_2$ & $\hat \alpha_3$ & $\hat \beta_1$ & $\hat \beta_2$ & $\hat \beta_3$ \\
\hline
Clenshaw-Curtis & $1.22$ & $1.53$ & $0.25$ & $-1.57$ & $-1.71$ & $ 0.25$ \\
Leja            & $1.18$ & $1.46$ & $0.35$ & $-2.37$ & $-2.40$ & $-0.76$ \\
\hline
\end{tabular}
\caption{Estimated $\hat{\bm \alpha}$ and $\hat{\bm \beta}$ parameters for $f_4(\bm y)$.}\label{table:burg}
\end{center}
\end{table}

%\newpage

\section{Conclusion}

In this work, we have presented an adaptive algorithm for constructing a sequence of interpolants of a function defined over a product of one dimensional intervals and admitting analytic extension to poly-ellipse in complex space. Following recent results in ``best $M$-term'' approximation, we derive a heuristic estimate for the optimal interpolation space, which we parametrize by two vectors, one depending on the region of analytic extension of the function and one depending on the Lebesgue constant of the interpolation scheme. Traditional methods for (quasi-)optimal approximation rely on a priori estimates, instead, we present a procedure for constructing a sequence of polynomial spaces, where each space is derived from an estimate inferred from an interpolant in the previous space. Each interpolant is constructed using a sparse grids approach, and we present a strategy for selecting the tensor rules so that the resulting grid is optimal (i.e., fewest nodes) with respect to the polynomial space. We also present several novel interpolation rules derived from greedy optimization of operator norms. Our numerical experiments demonstrate that the method can be applied to a wide range of problems without the need for a priori estimates of the region of analyticity of the function, so long as the function is analytic and the one dimensional family of rules inducing the sparse grid exhibits polynomial growth of the Lebesgue constant.

{\bf ACKNOWLEDGEMENTS}
This material is based upon work supported in part by the U.S. Air Force of Scientific Research under grant number 1854-V521-12 and by the U.S. Department of Energy, Office of Science, Office of Advanced Scientific Computing Research, Applied Mathematics program under contract and award numbers ERKJ259; and by the Laboratory Directed Research and Development program at the Oak Ridge National Laboratory, which is operated by UT-Battelle, LLC., for the U.S. Department of Energy under Contract DE-AC05-00OR22725.

\bibliographystyle{siam}
\bibliography{references}

\end{document}